\documentclass[leqno]{amsart}
\usepackage{cases}
\usepackage{cases}
\usepackage{amsmath}
\usepackage{graphicx}
\usepackage{mathrsfs}
\usepackage{bbm}
\usepackage{amsfonts}
\usepackage{amssymb}
\usepackage{amsmath,amsthm}
\usepackage{hyperref}
\usepackage{xcolor}
\usepackage{color}
\usepackage{ifpdf}
\usepackage{url}
\usepackage{enumerate}
\usepackage{multirow}
\usepackage{diagbox}
\usepackage{rotating}
\usepackage{tikz}

\definecolor{grey}{rgb}{0.5,0.5,0.4}
\def\opn#1#2{\def#1{\operatorname{#2}}} 
\opn\chara{char} \opn\length{\ell}
\opn\projdim{proj\,dim} \opn\injdim{inj\,dim} \opn\rank{rank}
\opn\depth{depth} \opn\grade{grade} \opn\height{height}
\opn\embdim{emb\,dim} \opn\codim{codim}

\opn\Tr{Tr} \opn\bigrank{big\,rank}
\opn\superheight{superheight}\opn\lcm{lcm}
\opn\trdeg{tr\,deg}%
\opn\reg{reg} \opn\lreg{lreg}
\opn\Ker{Ker} \opn\Coker{Coker} \opn\Im{Im} \opn\Hom{Hom}
\opn\Tor{Tor} \opn\Ext{Ext} \opn\End{End} \opn\Aut{Aut} \opn\id{id}

\opn\nat{nat}
\opn\pff{pf}
\opn\Pf{Pf} \opn\GL{GL} \opn\SL{SL} \opn\mod{mod} \opn\ord{ord}

\def\Implies{\ifmmode\Longrightarrow \else
     \unskip${}\Longrightarrow{}$\ignorespaces\fi}
\def\implies{\ifmmode\Rightarrow \else
     \unskip${}\Rightarrow{}$\ignorespaces\fi}
\def\iff{\ifmmode\Longleftrightarrow \else
     \unskip${}\Longleftrightarrow{}$\ignorespaces\fi}

\let\:=\colon

\newtheorem{Theorem}{Theorem}[section]
\newtheorem{Lemma}[Theorem]{Lemma}
\newtheorem{Corollary}[Theorem]{Corollary}

\newtheorem{Remark}[Theorem]{Remark}

\newtheorem{Definition}[Theorem]{Definition}
\newtheorem{Algorithm}[Theorem]{Algorithm}

\theoremstyle{definition}
\newtheorem*{pf}{Proof}
\opn\ini{in} \opn\inm{inm} \opn\Sym{Sym} \opn\diag{diag}
\opn\Ii{(i)} \opn\Iii{(ii)}

\title{ Convergence analysis of the Schwarz alternating method for unconstrained elliptic optimal control problems}
\author{Wei Gong $^\diamond$}
\author{Felix Kwok $^\dag$}
\author{Zhiyu Tan $^\ddag$}
\thanks{$^{\diamond}$ LSEC, Institute of Computational Mathematics, Academy of Mathematics and Systems
Science, Chinese Academy of Sciences, Beijing 100190, China.
Email: {\tt wgong@lsec.cc.ac.cn}. The author was supported by the Strategic Priority Research Program of Chinese Academy of Sciences (Grant No. XDB 41000000), the National Key Basic Research Program (Grant No. 2018YFB0704304) and the National Natural Science Foundation of China (Grant No. 12071468, 11671391).}
\thanks{$^{\dag}$ D\'{e}partment de math\'{e}matiques et de statistique, Universit\'{e} Laval, Qu\'{e}bec, Canada. Email: {\tt felix.kwok@mat.ulaval.ca} The author gratefully acknowledges support from the National Science and Engineering Research Council of Canada (RGPIN-2021-02595). The work described in this paper is partially supported by a grant from the ANR/RGC joint research scheme sponsored by the Research Grants Council of the Hong Kong Special Administrative Region, China and the French National Research Agency (Project no. A-HKBU203/19).}
\thanks{$^{\ddag}$ Center for Computation and Technology, Louisiana State University, Baton Rouge, LA 70803, USA. Email: {\tt ztan@cct.lsu.edu}}

\date{\today}

\begin{document}
\maketitle

{\bf Abstract:}\hspace*{10pt} {In this paper we analyze the Schwarz alternating method for unconstrained elliptic optimal control problems.  We discuss the convergence properties of the method in the continuous case first and then apply the arguments to the finite difference discretization case. In both cases, we prove that the Schwarz alternating method is convergent if its counterpart for an elliptic equation is convergent. Meanwhile, the convergence rate of the method for the elliptic equation under the maximum norm also gives a uniform upper bound (with respect to the regularization parameter $\alpha$) of the convergence rate of the method for the optimal control problem under the maximum norm of  proper error merit functions in the continuous case or vectors in the discrete case. Our numerical results corroborate our theoretical results and show that with $\alpha$ decreasing to zero, the method will converge faster. We also give some exposition of this phenomenon.}

{{\bf Keywords:}\hspace*{10pt}  The Schwarz alternating method,\ Saddle point problems,\ Elliptic optimal control problems,\ Maximum principle,\ Finite difference methods,\ Convergence rate.}
\vspace{8pt}


\section{Introduction}
\setcounter{equation}{0} 
In this paper, we consider the following
unconstrained elliptic distributed optimal  control problem
\begin{eqnarray}\label{OPTC}
\min\limits_{u\in L^2(\Omega)} \
J(y,u)=\frac{1}{2}\|y-y_d\|_{L^{2}(\Omega)}^2 + \frac{\alpha}{2}\|
u\|_{L^{2}(\Omega)}^2
\end{eqnarray}
subject to
\begin{equation}\label{OPT_PDE}
\left\{ \begin{aligned} \mathcal{L} y=f+u \ \ &\mbox{in}\
\Omega, \\
 \ y=0  \ \quad \quad &\mbox{on}\ \partial\Omega,\\
\end{aligned} \right.
\end{equation}
where $\Omega \subset \mathbb{R}^{d}\ (d = 1,2, 3)$ is a bounded Lipschitz domain,
$u\in L^{2}(\Omega)$ is the control variable, $y_d\in L^2(\Omega)$ is the
desired state or observation, $\alpha>0$ is the regularization
parameter and $\mathcal{L}$ is a self-adjoint and strictly elliptic operator which is defined as
\begin{equation}\label{operator_L}
\mathcal{L}y=-\sum\limits_{i,j=1}^{d}\frac{\partial }{\partial x_{j}}\left(a_{ij}(x)\frac{\partial y}{\partial x_{i}}\right) + c_{0}(x)y,
\end{equation}
with $a_{ij}(x),\ c_{0}(x)\in L^{\infty}(\Omega)$ and $c_{0}\geq 0$. 

The optimal control problem
(\ref{OPTC})-(\ref{OPT_PDE}) has a unique solution $u$ which can be
characterized by its first order optimality system (see, e.g.,
\cite{Lions_book}). Using standard arguments, the first order
necessary (also sufficient here) system of (\ref{OPTC})-(\ref{OPT_PDE}) is
\begin{equation}\label{COPTC}
\left\{
\begin{aligned}
&\alpha u+p=0 ;\\
 &\mathcal{L} y=f+u \ \ \ \mbox{in}\ \Omega, \quad y=0  \ \ \ \mbox{on}\ \partial\Omega;\\
&\mathcal{L} p=y-y_{d} \ \ \mbox{in}\ \Omega, \quad p=0  \ \ \ \mbox{on}\ \partial\Omega,\\
\end{aligned}
\right.
\end{equation}
where $p$ is the adjoint state.

By eliminating $u$, it is equivalent to 
\begin{equation}\label{COPTC1}
\left\{
\begin{aligned}
 &\mathcal{L}  y=f -\alpha^{-1}p \ \ \mbox{in}\ \Omega, \quad y=0  \ \ \ \mbox{on}\ \partial\Omega;\\
&\mathcal{L} p=y-y_{d} \ \ \ \ \ \ \mbox{in}\ \Omega, \quad p=0  \ \ \ \mbox{on}\ \partial\Omega,\\
\end{aligned}
\right.
\end{equation}
which can be treated as a saddle point problem.

In the past several decades, optimization problems with PDE constraint have attracted more and more attentions due to the increasingly broad applications in modern science and engineering.  The requirement of efficient solving of such kind of problems stimulates the research of this field as well. However, the fast and efficient simulations of such kind of problems are full of challenges.   For the theoretical results, we refer to \cite{Lions_book, Troltzsch_book}. A lot of researchers from different fields have also tried hard to deal with the related numerical problems, such as the optimization algorithms, the discretization methods and the fast solvers etc. and  fruitful achievements have been made. For the recent developments of the optimization algorithms and the convergence analysis of numerical schemes, we refer to the monograph \cite{ Hinze_book}.

In addition to the convergence analysis of the discretization schemes and the design of the optimization algorithms, the fast and robust solving of the resulting discretization problems is also very important. There are many attempts to deal with this kind of problems based on different approaches. 
In \cite{Biros_2005, Wei_2018, Schoberl_2007, Zulehner_2011}, the authors  used  preconditioned Krylov subspace methods to solve the first order optimality system by constructing some block preconditioners. In \cite{Borzi_2003, Borzi_2009, Brenner_2020, Schoberl_2011,  Simon_2009, Borzi_2008}, the authors used mutigrid methods to design fast solvers. Another strategy is to use domain decomposition methods to deal with optimal control problems (see e.g., \cite{Bartlett_2006, Benamou_1996, Benamou_1999, Ciaramella_2021, Chang_2011, Gander_2018, Heinkenschloss_2006, Heinkenschloss_2007, Nguyen_phd, Tan_2017, Tan_phd}). We also mention \cite{Prudencio_2006, Prudencio_2007, Yang_2013} for the parallel implementations of domain decomposition type algorithms. In \cite{Gander_2016, Kwok_2017, Gander_2020}, the authors discuss the time domain decomposition methods for time related optimal control problems.

In this paper, we focus on the Schwarz alternating method for the optimal control problem. The Schwarz alternating method was first introduced  by H. A. Schwarz in \cite{Schwarz_1870} to prove the existence and uniqueness of the solution of Laplace's equation on general domain with prescribed  Dirichlet boundary condition.  In the 1970s, it was used as a computing method to solve the partial differential equations, especially after P. L. Lions proving the convergence of the method based on the variational principle in \cite{PLions_1988} which simplifies the convergence analysis of the method. In \cite{PLions_1989}, P. L.  Lions also discussed the convergence properties of the method based on the maximum principle. For more about the origins and developments of the Schwarz alternating method, one can refer to \cite{Gander_1}. The generalization of the Schwarz alternating method in different directions yields a variety of domain decomposition methods. Domain decomposition methods (DDM for short) have been successfully used to construct fast solvers for the self-adjoint and positive definite partial differential equations. The essential parallel ability makes them attractive in applications. For more details on the design and the convergence analysis of DDM for the equations, we refer to the monograph \cite{Widlund_book}, the review papers \cite{Xu_1992, Xu_1998} and the references cited therein. As to the design of DDM for nonselfadjoint  or indefinite problems, we refer to \cite{Cai_1992, Cai_1993} and the references therein.

In contrast to the equation case where the DDM have fruitful theoretical results and satisfactory numerical performance,  in the optimal control problem case there are a few theoretical results of the DDM which are far from satisfactory.   But numerous numerical experiments in the literature have illustrated the efficiency and robustness of DDM for the optimal control problem. Roughly speaking, there are two categories of DDM for optimal control problems: the PDE-level DDM where DDM are only applied to the state equation and the adjoint state equation separately (cf. \cite{Chang_2011}) and the Optimization-level DDM where DDM are designed to decompose the optimal control problem directly to some optimization subproblems (cf. \cite{Benamou_1996, Benamou_1999, Gander_2018}). In \cite{Benamou_1996, Benamou_1999}, the authors propose a non-overlapping domain decomposition method based on the Robin type transmission conditions and prove the convergence of the algorithm with an absence of the convergence rate which is further discussed in \cite{Delourme_2019, Delourme_2021, Xu_2019} under the optimized Schwarz framework. In \cite{Gander_2018}, the authors prove the convergence of several non-overlapping domain decomposition methods and prove that for appropriate choice of the related parameters, the corresponding algorithms will be convergent in at most three iterations. To the best of our  knowledge, there are no robust theoretical convergence results of Optimization-level overlapping DDM for the optimal control problem in the literature. The existing theoretical analysis of the DDM for elliptic equations is hard to extend to the optimal control problem ( or PDE-constrained optimization) case because of the essential saddle point structure of the first order optimality system. It is even harder to obtain robust theoretical convergence results of the methods with respect to the regularization parameter $\alpha$. In this paper, we will prove the robust convergence results of the Schwarz alternating method for the elliptic optimal control problem by defining proper error merit functions (or vectors in the discrete case) which are related to the proper norms that are used in \cite{Brenner_2020, Wei_2018, Wei_2021, Tan_2017, Schoberl_2007, Zulehner_2011} and using the maximum principle of the elliptic operator. We will also give a uniform upper bound of the convergence rate. We hope that the results of this paper can give some insights of the DDM for the optimal control problem.

Our analysis is mainly based on the  weak maximum principle of second order elliptic operators. We state it in the following theorem and one can refer to \cite{Trudinger_book} for more details.
\begin{Theorem}\label{Th_weak_max}
\cite[Theorem 8.1]{Trudinger_book} Let $\Omega\subset\mathbb{R}^d\ (d=1,2,3)$ be a domain and $\phi\in H^1(\Omega)$ satisfy $\mathcal{L}\phi\leq 0\ (\geq 0)$ in $\Omega$. Then
\begin{equation}\nonumber
\sup\limits_{\Omega} \phi \leq \sup\limits_{\partial \Omega} \phi^+\quad \ (\inf\limits_{\Omega} \phi \geq \inf\limits_{\partial \Omega} \phi^-)
\end{equation}
where $\phi^+=\max\{\phi,0\}\quad (\phi^-=\min\{\phi,0\})$,  
\begin{equation}\nonumber
\sup\limits_{\partial \Omega} \phi= \inf\{c| \phi\leq c\ \mbox{on}\ \partial \Omega,\ c\in\mathbb{R}\}\  \mbox{and}\ \inf\limits_{\partial \Omega} \phi=-\sup\limits_{\partial \Omega} (-\phi).
\end{equation}
\end{Theorem}

The rest part of this paper is organized as follows. In section \ref{SAM_EOCP} we describe the Schwarz alternating method for the optimal control problem in detail. In section \ref{Con_C} we prove the convergence of the Schwarz alternating method and obtain the convergence rate under the maximum norm in the continuous case. We show that the convergence rate of the method for the optimal control problem is bounded by its counterpart in the elliptic equation case. We also obtain the convergence of the method  under the $L^2$ norm. Following similar discussions, in section \ref{Con_D} we show the same convergence properties of the method in the finite difference discretization case. In the last section we carry out some numerical experiments to illustrate our theoretical results. We also give more explanations in one dimensional case in the Appendix part. 

Throughout  this paper we will use the standard notations for differential operators, function spaces and norms that can be found in \cite{Adams_book, Trudinger_book}. We will denote $C$ a generic positive constant independent of $\alpha$ and the mesh size.

\section{The Schwarz alternating method for elliptic optimal control problems}\label{SAM_EOCP}
\setcounter{equation}{0}
In this section, we will extend the Schwarz alternating method for elliptic partial differential equations to the optimal control problem. The main issue here is how to properly define the subproblem on each subdomain that can guarantee the convergence of the method. 

\subsection{The subproblem} The subproblem that we will use is still an optimal control problem with PDE constraint. We will give the definition below. 

\begin{Definition}
Let $\omega\subset \mathbb{R}^d\ (d=1,2,3)$ be a bounded Lipschitz domain. $p_{\Gamma},\ y_{\Gamma}\in H^{1/2}(\partial \omega)$ are  given functions on $\partial \omega$.
The following optimal control problem is called a
Dirichlet-Dirichlet optimal control problem:
\begin{eqnarray}\label{DDOPC}
\min\limits_{u\in L^2(\omega)} \ J(y,u)=\frac{1}{2}\|
y-y_d\|_{L^{2}(\omega)}^2 + \frac{\alpha}{2}\|
u\|_{L^{2}(\omega)}^2-\int_{\partial \omega}\frac{\partial
y}{\partial n_{\mathcal{L}}}p_{\Gamma}d\sigma
\end{eqnarray}
subject to
\begin{equation}
\left\{ \begin{aligned} \mathcal{L} y=f+u \ \ &\mbox{in}\
\omega, \\
 \ y=y_{\Gamma}   \quad \quad &\mbox{on}\ \partial \omega,\\
\end{aligned} \right.
\end{equation}
where
\begin{equation}\label{eq:partial_nu}
\frac{\partial y}{\partial n_{\mathcal{L}}} =  \sum\limits_{i,j=1}^{d}a_{ij}(x)\frac{\partial y}{\partial x_{i}}\cos({\bf n},x_j),
\end{equation}
$\cos({\bf n},x_j)=j$-th direction cosine of ${\bf n}$ and ${\bf n}$ is the unit outward normal vector of $\partial \omega$.
\end{Definition}

A standard argument gives the first order necessary (also sufficient here) system of the Dirichlet-Dirichlet optimal control problem as
\begin{equation}\label{DDFC}
\left\{
\begin{aligned}
\alpha u+p&=0 \ \ \mbox{in}\ \omega;\\
\mathcal{L} y&=f+u \ \ \ \mbox{in}\ \omega,\  
 \ y=y_{\Gamma}  \ \ \ \mbox{on}\ \partial \omega;\\
\mathcal{L} p&=y-y_{d} \ \ \mbox{in}\
\omega, \ 
 \ p=p_{\Gamma}  \ \ \ \mbox{on}\ \partial \omega. \\
\end{aligned}
\right.
\end{equation}
Eliminating $u$ by the first equation, we have an equivalent system
\begin{equation}\label{DDFC_Sadddle}
\left\{
\begin{aligned}
\mathcal{L} y&=f - \alpha^{-1}p \ \ \ \mbox{in}\ \omega,\   
 \ y=y_{\Gamma}  \ \ \ \mbox{on}\ \partial \omega;\\
\mathcal{L} p&=y-y_{d} \ \ \ \ \ \ \ \mbox{in}\
\omega, \ 
 \ p=p_{\Gamma}  \ \ \ \mbox{on}\ \partial \omega. \\
\end{aligned}
\right.
\end{equation}
Note that it can also be reformulated as a saddle point problem.

\subsection{The Schwarz alternating method} 
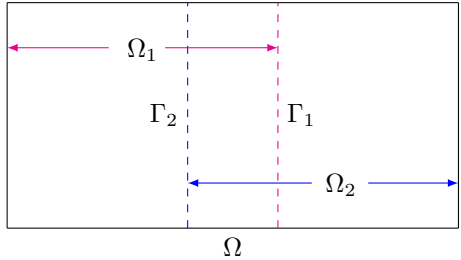
\begin{figure*}[htbp]
\begin{center}
\begin{tikzpicture}[scale = 3]
\draw (0,0) -- (2,0);
\draw (0,0) -- (0,1);
\draw (0,1) -- (2,1);
\draw (2,1) -- (2,0);
\node[below] at (1, 0) {$\Omega$};
\draw[dashed, color=blue] (0.8,0) -- (0.8,1);
\draw [-latex, color=blue] (1.6,0.2)--(2, 0.2);
\draw [-latex, color=blue] (1.35,0.2)--(0.8, 0.2);
\node[above] at (1.48,0.1) {$\Omega_2$};
\node[left] at (0.8, 0.5) {$\Gamma_2$};
\draw[dashed, color=magenta] (1.2,0) -- (1.2,1);
\draw [-latex, color=magenta] (0.7,0.8)--(1.2, 0.8);
\draw [-latex, color=magenta] (0.45,0.8)--(0, 0.8);
\node[above] at (0.6,0.7) {$\Omega_1$};
\node[right] at (1.2, 0.5) {$\Gamma_1$};
\end{tikzpicture}
\end{center}
\caption{Decomposition of $\Omega$}
\label{Fig:De_Omega}
\end{figure*}

Let $\{\Omega_i:\ i =1,2 \}$ be an overlapping domain decomposition  of $\Omega$ with $\Omega=\Omega_{1}\cup\Omega_{2}$ and $\Omega_{1}\cap\Omega_{2}\neq \emptyset$ (see Figure \ref{Fig:De_Omega}). Assume that $\Omega_i\ (i=1,2)$ are two bounded Lipschitz domains.

With these ingredients at hand, we can define the Schwarz alternating method for the elliptic optimal control problem. For the given overlapping decomposition $\Omega=\bigcup_{i=1}^{2}\Omega_{i}$, the Schwarz alternating method is described in Algorithm \ref{SAM}. We denote by $\frac{\partial y}{\partial n_{\mathcal{L},i}}\ (i=1,2)$ as that of $\frac{\partial y}{\partial n_{\mathcal{L}}}$ where $\omega$ is replaced by $\Omega_i\ (i=1,2)$.

For more about the Schwarz alternating method and domain decomposition methods of the optimal control problem, one can refer to \cite{Nguyen_phd, Tan_phd}.

\begin{Algorithm}\label{SAM}
(\textbf{The Schwarz Alternating Method})\\
1. Initialization: choose $y^{(0)},\ p^{(0)}\in H^1_0(\Omega)$. \\
2. For $k=0,1,\cdots$, solving the Dirichlet-Dirichlet optimal control problem on $\Omega_{i}\ (i=1,2)$ alternatively:\\
\begin{eqnarray}\nonumber
&\min\limits_{u^{(2k+i)}_{i}\in L^2(\Omega_{i})} \
J(y^{(2k+i)}_{i},u^{(2k+i)}_{i})=\frac{1}{2}\|
y^{(2k+i)}_{i}-y_d\|_{L^{2}(\Omega_{i})}^2 +
\frac{\alpha}{2}\|
u^{(2k+i)}_{i}\|_{L^{2}(\Omega_{i})}^2 \nonumber \\
&\quad\quad\quad\quad\quad\quad\quad\quad\quad\quad\quad-\int_{\partial \Omega_{i}}\frac{\partial
y^{(2k+i)}_{i}}{\partial n_{\mathcal{L},i}}[p^{(2k+i-1)}|_{\partial \Omega_{i}}]d\sigma\nonumber
\end{eqnarray}
\ \ \ \ subject to
\begin{equation}\nonumber
\left\{ \begin{aligned} \mathcal{L} y^{(2k+i)}_{i}=f+u^{(2k+i)}_{i} \ \
&\mbox{in}\
\Omega_{i}, \\
 \ y^{(2k+i)}_{i}=y^{(2k+i-1)}  \ \ \ \ &\mbox{on}\ \partial\Omega_{i}.\\
\end{aligned} \right.
\end{equation}
\\
3. Set $u^{(2k+i)}$, $y^{(2k+i)}$, $p^{(2k+i)}$ as follows:\\
\begin{equation}\nonumber
\begin{aligned}
u^{(2k+i)}& =\left\{
\begin{aligned}
&u^{(2k+i)}_{i}\ &\mbox{in}\ \Omega_{i},\\
&u^{(2k+i-1)}\ &\mbox{in}\ \overline{\Omega}\setminus\Omega_{i};\\
\end{aligned}
\right.\\
y^{(2k+i)}& = \left\{
\begin{aligned}
&y^{(2k+i)}_{i}\ &\mbox{in}\ \Omega_{i},\\
&y^{(2k+i-1)}\ &\mbox{in}\ \overline{\Omega}\setminus\Omega_{i};\\
\end{aligned}
\right.\\
p^{(2k+i)}&=-\alpha u^{(2k+i)}. \\
\end{aligned}
\end{equation}
\end{Algorithm}

\begin{Remark}
Since each problem on the subdomain is still an optimal control problem, the decomposition here is an Optimization-level decomposition of the original problem. Note that it decomposes the state  equation and the objective function simultaneously, which is one another feature of this algorithm.
\end{Remark}

\subsection{An equivalent algorithm for a saddle point problem} By the equivalence of the optimal control problems and their first order optimality systems, we can give the equivalent Schwarz alternating method for the first order optimality system (\ref{COPTC1}) in a parallel way. We will give it in Algorithm \ref{SAM_2}.

\begin{Algorithm}\label{SAM_2}
\ \\
1. Initialization: choose $y^{(0)},\ p^{(0)}\in H^1_0(\Omega)$.\\
2. For $k=0,1,\cdots$, solving the following problem on $\Omega_{i}\ (i=1,2)$ alternatively:\\
\begin{equation}\label{sub_system}
\left\{\begin{array}{llr}
\mathcal{L} y^{(2k+i)} = f-\alpha^{-1} p^{(2k+i)}  \quad&\mbox{in}\
\Omega_{i}, \\
\mathcal{L} p^{(2k+i)} = y^{(2k+i)} - y_{d} \quad&\mbox{in}\
\Omega_{i}, \\
\quad \ \ y^{(2k+i)}=y^{(2k+i-1)} \quad &\mbox{on}\ \partial \Omega_{i},\\
\quad \ \ p^{(2k+i)}=p^{(2k+i-1)} \quad &\mbox{on}\ \partial \Omega_{i},\\
\quad \ \ y^{(2k+i)}=y^{(2k+i-1)}  \quad &\mbox{in}\ \overline{\Omega}\setminus \overline{\Omega_{i}},\\ 
\quad \ \ p^{(2k+i)}=p^{(2k+i-1)} \quad &\mbox{in}\ \overline{\Omega}\setminus \overline{\Omega_{i}}.\\ 
\end{array}\right.  
\end{equation}
\end{Algorithm}

\begin{Remark}
The algorithm here can be viewed as a direct decomposition of a saddle point problem into two saddle points subproblems.
\end{Remark}

\section{Convergence analysis of the algorithm in the continuous case}\label{Con_C}
\setcounter{equation}{0}
In this section, we focus on the convergence analysis of Algorithm \ref{SAM_2}, while the convergence properties of Algorithm \ref{SAM} follow directly from the equivalence of these two algorithms. Our purpose is to prove the robust convergence properties of the algorithm with respect to the parameter $\alpha$.

Before we give the analysis in detail, we first give a key observation of the self-adjoint and strictly elliptic operator $\mathcal{L}$.
\begin{Lemma}\label{lem:L_pro}
For a bounded domain $\omega\subset \mathbb{R}^d\ (d=1,2,3)$ and $\beta>0$, we assume $\phi,\psi\in H^1(\omega)$ satisfy 
\begin{equation}\nonumber
\left\{
\begin{aligned}
&\mathcal{L}\psi = - \beta\phi\ \ \mbox{in}\ \ \omega;\\
&\mathcal{L}\phi = \psi\ \ \ \ \ \ \mbox{in}\ \ \omega.
\end{aligned}
\right.
\end{equation}
Then we have
\begin{equation}\nonumber
\mathcal{L}(\psi^2 + \beta \phi^2)  \leq 0\quad \mbox{in}\ \ \omega.
\end{equation}
\end{Lemma} 
\begin{pf}
By the definition of $\mathcal{L}$, we have
\begin{equation}\nonumber
\begin{aligned}
\mathcal{L}(\psi^2) 
&= -\sum\limits_{i,j=1}^{d}\frac{\partial }{\partial x_{j}}\left(a_{ij}(x)\frac{\partial (\psi^2)}{\partial x_{i}}\right) + c_{0}(x)\psi^2\\
& = -2\psi\sum\limits_{i,j=1}^{d}\frac{\partial }{\partial x_{j}}\left(a_{ij}(x)\frac{\partial \psi}{\partial x_{i}}\right) - 2\sum\limits_{i,j=1}^{d}\left(a_{ij}(x)\frac{\partial\psi}{\partial x_{i}}\frac{\partial \psi}{\partial x_{j}}\right)+ c_{0}(x)\psi^2\\
& = 2\psi\mathcal{L}\psi - 2\sum\limits_{i,j=1}^{d}\left(a_{ij}(x)\frac{\partial\psi}{\partial x_{i}}\frac{\partial \psi}{\partial x_{j}}\right) - c_{0}(x)\psi^2.
\end{aligned}
\end{equation}
Therefore, the strictly ellipticity of $\mathcal{L}$ implies
\begin{equation}\nonumber
\begin{aligned}
\mathcal{L}(\psi^2 + \beta\phi^2) & = 2\psi\mathcal{L}\psi  - 2\sum\limits_{i,j=1}^{d}\left(a_{ij}(x)\frac{\partial\psi}{\partial x_{i}}\frac{\partial \psi}{\partial x_{j}}\right) - c_{0}(x)\psi^2 \\
&\quad + \beta\left[ 2\phi\mathcal{L}\phi  - 2\sum\limits_{i,j=1}^{d}\left(a_{ij}(x)\frac{\partial\phi}{\partial x_{i}}\frac{\partial \phi}{\partial x_{j}}\right) - c_{0}(x)\phi^2\right]\\
& = - 2\sum\limits_{i,j=1}^{d}\left(a_{ij}(x)\frac{\partial\psi}{\partial x_{i}}\frac{\partial \psi}{\partial x_{j}}\right)  - 2\beta\sum\limits_{i,j=1}^{d}\left(a_{ij}(x)\frac{\partial\phi}{\partial x_{i}}\frac{\partial \phi}{\partial x_{j}}\right) - c_{0}(x)(\psi^2 + \beta \phi^2)\\
&\leq 0.
\end{aligned}
\end{equation}
This finishes the proof of the lemma.
\end{pf}

Now we prove the convergence of Algorithm \ref{SAM_2} under the maximum norm by applying the weak maximum principle of Theorem \ref{Th_weak_max}. We first give the error equations and then introduce some auxiliary problems. By investigating the relationship between the errors and the solutions of the auxiliary problems, we will give the convergence of the algorithm.

Denote $e_{y}^{(j)}=y-y^{(j)}$ and $e_p^{(j)}=p-p^{(j)}$ with $y,\ p$ the solutions of (\ref{COPTC1}) and $y^{(j)},\ p^{(j)}$ generated by Algorithm \ref{SAM_2}. A direct calculation shows that the errors $e_{y}^{(j)}$, $e_{p}^{(j)}\ ( j=2k+i, i=1,2, k=0,1,2,...$) satisfy the following equations
\begin{equation}\label{error_system}
\left\{\begin{array}{llr}
\mathcal{L} e_{y}^{(2k+i)}  = -\alpha^{-1} e_{p}^{(2k+i)} \quad&\mbox{in}\
\Omega_{i}, \\
\mathcal{L} e_{p}^{(2k+i)} = e_{y}^{(2k+i)} \quad &\mbox{in}\
\Omega_{i}, \\
\ \ e_{y}^{(2k+i)} =e_{y}^{(2k+i-1)} \quad &\mbox{on}\ \partial \Omega_{i},\\
\ \ e_{p}^{(2k+i)} =e_{p}^{(2k+i-1)} \quad &\mbox{on}\ \partial \Omega_{i},\\
\ \ e_{y}^{(2k+i)} =e_{y}^{(2k+i-1)}  \quad &\mbox{in}\ \overline{\Omega}\setminus \overline{\Omega_{i}},\\ 
\ \  e_{p}^{(2k+i)} =e_{p}^{(2k+i-1)} \quad &\mbox{in}\ \overline{\Omega}\setminus \overline{\Omega_{i}}.\\ 
\end{array}\right.  
\end{equation}

We denote
\begin{equation}\label{eta}
\eta^{(j)}=(e_y^{(j)} )^2 + \alpha^{-1} (e_p^{(j)})^2\geq 0.
\end{equation}
with $j=2k+i\ ( i=1,2, k=0,1,...$). For each $j=2k+i$, by Lemma \ref{lem:L_pro}, $\eta^{(j)}$ satisfies 
\begin{equation}\label{eta_system}
\left\{\begin{array}{llr}
\quad \mathcal{L} \eta^{(2k+i)} \leq 0 \quad&\mbox{in}\
\Omega_{i}, \\
\quad \ \ \eta^{(2k+i)}=\eta^{(2k+i-1)} \quad &\mbox{on}\ \partial \Omega_{i},\\
\quad \ \ \eta^{(2k+i)}=\eta^{(2k+i-1)}  \quad &\mbox{in}\ \overline{\Omega}\setminus \overline{\Omega_{i}}.\\ 
\end{array}\right.  
\end{equation}

In the following, we introduce some auxiliary problems. We set $\xi^{(j)}$ satisfies 
\begin{equation}\label{xi_system1}
\left\{\begin{array}{llr}
\quad \mathcal{L} \xi^{(2k+i)} = 0 \quad&\mbox{in}\
\Omega_{i}, \\
\quad \ \ \xi^{(2k+i)}=\xi^{(2k+i-1)} \quad &\mbox{on}\ \partial \Omega_{i},\\
\quad \ \ \xi^{(2k+i)}=\xi^{(2k+i-1)}  \quad &\mbox{in}\ \overline{\Omega}\setminus \overline{\Omega_{i}}\\ 
\end{array}\right.  
\end{equation}
with $j = 2k + i$, $i=1,2$, $k = 0,1,2,\dots$, $\xi^{(0)}=\eta^{(0)}$ on $\overline{\Omega}$. 

\begin{Remark}\label{rem:SAM_eq}
It is worth to noting that the system (\ref{xi_system1}) is exactly the $k$-th iteration of the  Schwarz alternating method for the elliptic equation 
$$\mathcal{L}\xi=0\ \ \mbox{in}\ \ \Omega\quad \mbox{and} \quad \xi=0\ \ \mbox{on}\ \ \partial \Omega$$
on $\Omega_{i}$ with the initial guess $\xi^{(0)}=\eta^{(0)}$. 
\end{Remark}

\begin{Lemma}\label{lem:E_OCP_Eq}
For $ i = 1,2$ and $k = 0,1,2,\dots$, suppose that $\eta^{(2k+i)}$ and $\xi^{(2k+i)}$ are defined as above, then
\begin{equation}\nonumber
0\leq \eta^{(2k+i)} \leq \xi^{(2k+i)}\ \ \mbox{in} \ \overline{\Omega_{i}}.
\end{equation}
\end{Lemma}

\begin{pf}
We will prove it by induction.

According to (\ref{eta_system}), (\ref{xi_system1}) and $\xi^{(0)}=\eta^{(0)}$ on $\overline{\Omega}$, we have
\begin{equation}\nonumber
\left\{\begin{array}{llr}
 \mathcal{L} \left(\eta^{(1)}-\xi^{(1)}\right) \leq 0 \quad&\mbox{in}\
\Omega_{i}, \\
\quad\ \ \ \eta^{(1)}-\xi^{(1)} = 0\quad &\mbox{on}\ \partial \Omega_{i},\\
\end{array}\right.  
\end{equation}
which follows by Theorem \ref{Th_weak_max} that
$$ 0\leq \eta^{(1)} \leq \xi^{(1)}\ \ \mbox{in} \ \overline{\Omega_{1}},$$
i.e., the results hold for $k = 0$ and $i=1$.

Without loss of generality, we assume the results hold for $2k+i-1$, then we prove that the results hold for $2k+i$. 

Since
$\eta^{(2k+i)}- \xi^{(2k+i)}$ satisfies the following equation
\begin{equation}\nonumber
\left\{\begin{array}{llr}
 \mathcal{L} \left(\eta^{(2k+i)}-\xi^{(2k+i)}\right) \leq 0 \quad&\mbox{in}\
\Omega_{i}, \\
\quad\ \ \ \eta^{(2k+i)}-\xi^{(2k+i)}=\eta^{(2k+i-1)} -\xi^{(2k+i-1)}\quad &\mbox{on}\ \partial \Omega_{i}.\\
\end{array}\right.  
\end{equation}
According to the assumption of the induction, we know that 
$$\eta^{(2k+i-1)} -\xi^{(2k+i-1)}\leq 0.$$ 
By Theorem \ref{Th_weak_max}, we have 
\begin{equation}\nonumber
\eta^{(2k+i)}-\xi^{(2k+i)}\leq 0\ \ \mbox{in}\ \ \overline{\Omega_{i}},\quad \mbox{i.e.,}\quad  \eta^{(2k+i)} \leq \xi^{(2k+i)}\ \ \mbox{in}\ \ \overline{\Omega_{i}}.
\end{equation}
Now, we complete the proof of the lemma.
\end{pf}

Recalling Remark \ref{rem:SAM_eq}, Lemma \ref{lem:E_OCP_Eq} gives the relationship of the errors $\eta^{(j)}$ of Algorithm \ref{SAM_2} and the solutions $\xi^{(j)}$, i.e.,
\begin{equation}\nonumber
0\leq\eta^{(j)}\leq \xi^{(j)}\ \ \mbox{on}\ \ \overline{\Omega},\quad \forall\ j = 2k+i, k= 0,1,2,\cdots, i = 1,2,
\end{equation}
if $\xi^{(0)} = \eta^{(0)}$ on $ \overline{\Omega}$. This will give the convergence of Algorithm \ref{SAM_2} if $\xi^{(j)}\rightarrow 0$ as $j\rightarrow +\infty$. Furthermore, we can also obtain a uniform upper bound of the convergence rate of Algorithm \ref{SAM_2}. We will prove it and collect all these convergence results in the following main theorem of this paper.

\begin{Theorem}\label{Th:MN_C}
For a given domain decomposition, if the Schwarz alternating method for 
\begin{equation}\label{eq:L}
\mathcal{L} \xi=0\ \ \mbox{in}\ \ \Omega\quad \mbox{and}\quad \xi=0\ \ \mbox{on}\ \ \partial \Omega
\end{equation}
is convergent, then the Schwarz alternating algorithm, i.e., Algorithm \ref{SAM_2} for the system (\ref{COPTC1}) is also convergent. Moreover, if the convergence rate of the Schwarz alternating method for the elliptic equation (\ref{eq:L}) is $\rho_e\in(0,1)$ under the maximum norm, i.e.,
\begin{equation}\label{eq:SAMEq_rate}
\sup\limits_{x\in \overline{\Omega}} \xi^{(2k)}\leq \rho_{e} \sup\limits_{x\in \overline{\Omega}}\xi^{(2(k-1))}\quad k = 1,2,\cdots,
\end{equation}
then for $k=1, 2, ...$,
we have
\begin{equation}\label{eq:SAM_rate}
\sup\limits_{x\in \overline{\Omega}} \eta^{(2k)}\leq \rho_{e} \sup\limits_{x\in \overline{\Omega}}\eta^{(2(k-1))}.
\end{equation}
\end{Theorem}
\begin{pf}
We only need to prove (\ref{eq:SAM_rate}).

For $k\geq 1$, we take the $2(k-1)$-th step as the initial step, i.e.,
$$\xi^{(2(k-1))} =  \eta^{(2(k-1))}\ \ \mbox{on}\ \ \overline{\Omega}$$
in (\ref{eta_system}) and (\ref{xi_system1}).
Applying Lemma \ref{lem:E_OCP_Eq}, we have
$$ \eta^{(2k)} \leq \xi^{(2k)}\ \ \mbox{on}\ \ \overline{\Omega}.$$


Therefore, combining with (\ref{eq:SAMEq_rate}), we have
$$ \sup\limits_{x\in \overline{\Omega}} \eta^{(2k)}\leq \sup\limits_{x\in \overline{\Omega}} \xi^{(2k)}\leq \rho_e \sup\limits_{x\in \overline{\Omega}}\xi^{(2(k-1))} = \rho_e \sup\limits_{x\in \overline{\Omega}}\eta^{(2(k-1))},$$
which completes the proof.
\end{pf}

By applying the above theorem, we can also obtain the convergence of the algorithm under the $L^2$ norm.
\begin{Corollary} \label{cor:L2_C}
Suppose the assumptions in Theorem \ref{Th:MN_C} hold, then we have 
\begin{equation}\nonumber
\|e_{y}^{(2k)}\|_{0}^2+ \alpha^{-1}\|e_{p}^{(2k)}\|_{0}^2\leq C\rho_e^k\sup\limits_{x\in \overline{\Omega}} \eta^{(0)},
\end{equation}
where $C>0$ is a constant independent of $\rho_e$ and $\alpha$.
\end{Corollary}
\begin{pf}
According to Theorem \ref{Th:MN_C}, we have
\begin{equation}\nonumber
\|e_{y}^{(2k)}\|_{0}^2+\alpha^{-1}\|e_{p}^{(2k)}\|_{0}^2 = \int_{\Omega} \eta^{(2k)}dx \leq |\Omega| \sup\limits_{x\in \overline{\Omega}} \eta^{(2k)}\leq C \rho_e^k \sup\limits_{x\in \overline{\Omega}} \eta^{(0)}.
\end{equation}
\end{pf}

\begin{Remark}\label{rem:Th_rate}
A number of remarks are given as follows:
\begin{enumerate}
\item One can refer to \cite{Schwarz_1870, PLions_1989} for the convergence analysis and the convergence rate of the Schwarz alternating method for the elliptic equations under the maximum norm. 

\item Estimate (\ref{eq:SAM_rate}) gives a uniform upper bound, i.e.,  $\rho_e$,  for the convergence rate of the method for the optimal control problem. Our numerical results in Section 5 and the theoretical analysis for the 1$D$ continuous case in the appendix part indicate this upper bound is not optimal and the optimal one should be $\rho_e^2$.


\item The choice of $\alpha$ will affect the convergence rate of Algorithm \ref{SAM_2}. Roughly speaking, as $\alpha$ becomes smaller, the algorithm will converge faster. Our numerical results in section \ref{NE} illustrate this. In the one dimensional case, we will attach more theoretical results in the appendix part of this paper to explain this. One can also refer to \cite{Tan_phd} for more details.

\item The convergence rate of the Schwarz alternating method for the elliptic equation
\begin{equation}\nonumber
\mathcal{L}\tilde{\xi} + 2\alpha^{-\frac{1}{2}}\tilde{\xi} = 0\ \ \mbox{in}\ \ \Omega\quad \mbox{and}\quad \xi=0\ \ \mbox{on}\ \ \partial \Omega
\end{equation}
is a much better estimate of the convergence rate of the method for the optimal control problem for $\alpha$ small. Our numerical results and the discussion in $1$D case (see Appendix B) corroborate this claim.

\end{enumerate}

\end{Remark}

\section{Extension to the discrete case}\label{Con_D}
\setcounter{equation}{0}
In this section, we will consider the possibility to extend the arguments in the previous section to the discrete case. 

For an illustration purpose, we will consider the problem in two dimensional case and use finite difference methods to discretize the optimal control problem. More precisely, we will take 
$$\Omega = (0,1)^2\subset\mathbb{R}^2 \quad \mbox{and} \quad\mathcal{L} = -\Delta.$$

We first give the discrete problem of the optimal control problem discretized by the finite difference method and then give the Schwarz alternating method for the discrete problem. After that, we will prove the convergence results of the method following the arguments in parallel with that of its continuous counterpart.

\subsection{The discrete problem}
We adopt the five-point finite difference scheme to discretize the state equation and approximate the objective functional  by the trapezoidal rule. We refer to \cite{Borzi_2003, Liu_2019} for the numerical analysis of the finite difference scheme for the optimal control problem. 

For a given step size $h=1/N$, $\Omega=(0,1)^2$ is partitioned uniformly by grid $(x_{i}, y_{j})$ with $0\leq i, j\leq N$ and $x_{i}=ih,\ y_{j}=jh$. For given function $z$, we take $z_{i,j} = z(x_i,y_j)$ or some approximation of it and set $\{z_{i,j}:\ 0\leq i,j\leq N\}$, which can be mapped to a vector $Z\in \mathbb{R}^{(N+1)^2}$. We denote by $Z_{I}$ the part corresponding to the interior grid points in $\Omega$ and $Z_{\partial}$ the part corresponding to the boundary grid points on $\partial \Omega$.

The five-point finite difference scheme for the state equation is given by
\begin{equation}\label{FDM_1}
-\frac{y_{i-1,j}+y_{i+1,j}-4y_{i,j}+y_{i,j-1}+y_{i,j+1}}{h^2}=f_{i,j} + u_{i,j}\quad 1\leq i,j \leq N-1,
\end{equation} 
and the discrete problem of (\ref{OPTC})-(\ref{OPT_PDE}) reads
\begin{eqnarray}\label{dis_OPTC}
\min\limits_{U\in \mathbb{R}^n} \
J(Y,U)=\frac{h^2}{2}\|Y_{I}-Y_{d,I}\|_{\mathbb{R}^n}^2 + \frac{\alpha h^2}{2}\|
U_{I}\|_{\mathbb{R}^n}^2
\end{eqnarray}
subject to
\begin{equation}\label{dis_OPT_PDE}
\left\{
\begin{aligned}
\mathcal{L}_{h}Y&=F_{I}+U_{I},\\
Y_{\partial}&=0,
\end{aligned}
\right.
\end{equation}
where $\mathcal{L}_{h}$ is the corresponding coefficient matrix of the five-point finite difference scheme and $n=(N-1)^2$.

Following a standard argument, the first order optimality system of the discrete problem  reads
\begin{equation}\label{dis_FOC}
\left\{
\begin{aligned}
& \mathcal{L}_{h}Y=F_{I}+U_{I},\ \ Y_{\partial}=0;\\
& \mathcal{L}_{h}P=h^2(Y_{I}-Y_{d,I}),\ \ P_{\partial}=0;\\
&\alpha h^2 U_{I} + P_{I} =0.
\end{aligned}
\right.
\end{equation}
Similarly, we can eliminate $U_{I}$ and after doing some reformulation, we have
\begin{equation}\label{dis_FOC2}
\left\{
\begin{aligned}
& \mathcal{L}_{h}Y=F_{I}- \frac{1}{\alpha }(\frac{1}{h^2}P_{I}),\ \ Y_{\partial}=0;\\
& \mathcal{L}_{h}(\frac{1}{h^2}P)=Y_{I}-Y_{d,I},\ \ P_{\partial}=0.
\end{aligned}
\right.
\end{equation}

\begin{Remark}\label{Remark_1}
We use the discretize-then-optimize approach in the above setting. If we use the optimize-then-discretize approach, we can obtain the following linear system
 \begin{equation}\label{dis_od_FOC}
\left\{
\begin{aligned}
& \mathcal{L}_{h}Y=F_{I}+U_{I},\ \ Y_{\partial}=0;\\
& \mathcal{L}_{h}\tilde{P}=Y_{I}-Y_{d,I},\ \ \tilde{P}_{\partial}=0;\\
&\alpha U_{I} + \tilde{P}_{I} =0,
\end{aligned}
\right.
\end{equation}
which is the first order optimality system of the following discrete optimization problem
\begin{eqnarray}\label{dis_OPTC2}
\min\limits_{U\in \mathbb{R}^n} \
J(Y,U)=\frac{1}{2}\|Y_{I}-Y_{d,I}\|_{\mathbb{R}^n}^2 + \frac{\alpha}{2}\|
U_{I}\|_{\mathbb{R}^n}^2
\end{eqnarray}
subject to
\begin{equation}\label{dis_OPT_PDE2}
\left\{
\begin{aligned}
\mathcal{L}_{h}Y&=F_{I}+U_{I},\\
Y_{\partial}&=0.
\end{aligned}
\right.
\end{equation}

Eliminating $U_{I}$, we have
\begin{equation}\label{dis_od_FOC2}
\left\{
\begin{aligned}
& \mathcal{L}_{h}Y=F_{I}- \frac{1}{\alpha }\tilde{P}_{I},\ \ Y_{\partial}=0;\\
& \mathcal{L}_{h}\tilde{P}=Y_{I}-Y_{d,I},\ \ \tilde{P}_{\partial}=0.\\
\end{aligned}
\right.
\end{equation}
Note that $\tilde{P}=\frac{1}{h^2}P$. More details about these schemes can be found in \cite{Liu_2019}.
\end{Remark}

\subsection{The Schwarz alternating method for the discrete problem}
As in the continuous case, we can define the Schwarz alternating algorithms analogous to Algorithm \ref{SAM} and Algorithm \ref{SAM_2} in this discrete setting. We give the counterpart of Algorithm \ref{SAM_2} and prove the convergence of the algorithm. The convergence properties of the equivalent algorithm in the discrete case, which corresponds to Algorithm \ref{SAM},  follow directly.

For $i=1,2$, we denote by $Z_{I},\ Z_{\partial},\ Z_{i}, \ Z_{i,I},\ Z_{i,\partial}$ the vector components of $Z$ corresponding to the grid points related to $\Omega,\ \partial \Omega,\ \overline{\Omega_{i}},\ \Omega_{i},\ \partial\Omega_{i}$, respectively, and $\mathcal{L}_{h}^{(i)}$ the restriction of $\mathcal{L}_{h}$ to the grid points related  to $\overline{\Omega_{i}}$.   

\begin{Algorithm}\label{dis_SAM_2}
\ \\
1. Initialization: choose $Y^{(0)},\ P^{(0)}\in \mathbb{R}^{(N+1)^2}$ with $Y^{(0)}_{\partial}=0$ and $P^{(0)}_{\partial}=0$.\\
2. For $k=0,1,\cdots$, solving the following problem alternatively:\\
\begin{equation}\label{dis_sub_system}
\left\{\begin{aligned}
\mathcal{L}_{h}^{(i)} Y^{(2k+i)}_{i} &= F_{i,I}-\alpha^{-1}(\frac{1}{h^2} P^{(2k+i)}_{i,I}), \\
\mathcal{L}_{h}^{(i)}( \frac{1}{h^2} P^{(2k+i)}_{i} ) &= Y_{i,I}^{(2k+i)} - Y_{d,i,I}, \\
\quad \ \ Y^{(2k+i)}_{i,\partial} &=Y^{(2k+i-1)}_{i,\partial},\\
\quad \ \ P^{(2k+i)}_{i,\partial} &=P^{(2k+i-1)}_{i,\partial} ,\\
\quad \ \ Y^{(2k+i)} &=Y^{(2k+i-1)}\ \mbox{others},\\ 
\quad \ \ P^{(2k+i)} &=P^{(2k+i-1)}\ \mbox{others}.
\end{aligned}\right.  
\end{equation}
\end{Algorithm}
\begin{Remark}
On each subdomain, the discrete subproblem in the above algorithm is the first order optimality system of the discrete problem of the Dirichlet-Dirichlet optimal control problem on it.
\end{Remark}

\subsection{Convergence analysis}
In the following, we prove the convergence results of Algorithm \ref{dis_SAM_2}. We first write down the error equations of the algorithm and then prove a discrete counterpart of Lemma \ref{lem:L_pro} for $\mathcal{L}_h$ for a proper error merit vector. The remaining analysis can be carried out identically in parallel with the continuous case by recalling the maximum principle of $\mathcal{L}_h$. 

\subsubsection{The error equations}
We define $E_{y}^{(j)}=Y-Y^{(j)}$, $E_{p}^{(j)}=P-P^{(j)}$ ($j = 2k+i,\ i=1,2,\ k=0,1,2,\dots$) with $Y,\ P$ the solutions of (\ref{dis_FOC2}) and $Y^{(j)},\ P^{(j)}$ generated by Algorithm \ref{dis_SAM_2}. Then the errors $E_{y}^{(j)}$ and $E_{p}^{(j)}$ satisfy 
\begin{equation}\label{dis_error_sub}
\left\{\begin{aligned}
\mathcal{L}_{h}^{(i)} E^{(2k+i)}_{y,i} &= -\alpha^{-1}(\frac{1}{h^2} E^{(2k+i)}_{p,i,I}), \\
\mathcal{L}_{h}^{(i)}( \frac{1}{h^2} E^{(2k+i)}_{p,i} ) &= E_{y,i,I}^{(2k+i)} , \\
\quad \ \ E^{(2k+i)}_{y,i,\partial} &=E^{(2k+i-1)}_{y,i,\partial},\\
\quad \ \ E^{(2k+i)}_{p,i,\partial} &=E^{(2k+i-1)}_{p,i,\partial} ,\\
\quad \ \ E^{(2k+i)}_{y} &=E^{(2k+i-1)}_{y}\ \mbox{others},\\ 
\quad \ \ E^{(2k+i)}_{p} &=E^{(2k+i-1)}_{p}\ \mbox{others}.\\ 
\end{aligned}\right.  
\end{equation}

\subsubsection{A property of $\mathcal{L}_h$}
 For $Z=(Z_{k})_{m\times 1}\in \mathbb{R}^m$, we define $Z^2\in \mathbb{R}^m$ in a componentwise sense, i.e., $(Z^2)_{k}=(Z_{k})^2$. 
 
\begin{Lemma}\label{lem:Lh_pro}
Let $\beta>0$ and $\Phi,\Psi\in \mathbb{R}^{(N+1)^2}$ satisfy
\begin{equation}\nonumber
\left\{
\begin{aligned}
&\mathcal{L}_h\Psi = -\beta \Phi_I,\\
&\mathcal{L}_h\Phi = \Psi_I,
\end{aligned}
\right.
\end{equation}
then
\begin{equation}\nonumber
\mathcal{L}_h(\Psi^2 + \beta \Phi^2)\leq 0,
\end{equation}
which is understood componentwisely.
\end{Lemma} 
\begin{pf}
Let $\{\phi_{i,j}:\ 0\leq i, j\leq N \}$ be the set corresponding to $\Phi$. Since
\begin{equation}\nonumber
\begin{aligned}
&\quad-\frac{\phi_{i-1,j}^2+\phi_{i+1,j}^2-4\phi_{i,j}^2+\phi_{i,j-1}^2+\phi_{i,j+1}^2}{h^2}\\
&=-2\phi_{i,j}\frac{\phi_{i-1,j}+\phi_{i+1,j}-4\phi_{i,j}+\phi_{i,j-1}+\phi_{i,j+1}}{h^2}\\
&\ \ \ \ \ -\frac{(\phi_{i-1,j}-\phi_{i,j})^2+(\phi_{i+1,j}-\phi_{i,j})^2}{h^2}-\frac{(\phi_{i,j-1}-\phi_{i,j})^2+(\phi_{i,j+1}-\phi_{i,j})^2}{h^2},
\end{aligned}
\end{equation}
we have
$$ \mathcal{L}_h(\Phi^2) \leq -2\Phi_I^T\mathcal{L}_h\Phi.$$
It follows
$$ \mathcal{L}_h(\Psi^2 + \beta \Phi^2) \leq -2\Psi_I^T\mathcal{L}_h\Psi - 2\beta\Phi_I^T\mathcal{L}_h\Phi=0$$
by the assumption. This completes the proof.
\end{pf}

\begin{Remark}\label{rem:LtoLi}
The above lemma also applies to  $\mathcal{L}_h^{(i)}\ (i=1,2)$ for vectors related to $\overline{\Omega_i}$.
\end{Remark}

\subsubsection{The maximum principle for $\mathcal{L}_h$} 
The maximum principle plays a crucial role in the proof of the convergence properties of the Schwarz alternating algorithm in the continuous case. For the finite difference operator $\mathcal{L}_{h}$, we have a similar discrete maximum principle. We state it in the following theorem and refer to \cite[Theorem 3]{Ciarlet_1970} for more details.

\begin{Theorem}\label{dis_weak_max}
The finite difference operator $\mathcal{L}_{h}$ satisfies the discrete maximum principle, i.e., for any $Z\in \mathbb{R}^{(N+1)^2}$ with $\mathcal{L}_{h}Z\leq 0\ (\geq 0)$, we have
\begin{equation}\nonumber
\max\{Z\}\leq \max\{0,Z_{\partial }\}\quad (\min\{Z\}\geq \min\{0,Z_{\partial }\}),
\end{equation}
where $\max$ and $\min$ are understood in a componentwise sense.
\end{Theorem}

\begin{Remark}\label{rem:LtoLi}
Note that this theorem also applies to  $\mathcal{L}_h^{(i)}\ (i=1,2)$ for vectors related to $\overline{\Omega_i}$.
\end{Remark}

\subsubsection{The convergence results}
We define the error merit vectors as
\begin{equation}\nonumber
E^{(2k+i)}=(E_{y}^{(2k+i)})^2+\frac{\alpha^{-1}}{h^4}(E_{p}^{(2k+i)})^2,\ \ i=1,2,\ k = 0,1,2,\dots.
\end{equation}
Applying Lemma \ref{lem:Lh_pro} to $\mathcal{L}_{h}^{(i)}$ and $E^{(2k+i)}_{i}$, we have
\begin{equation}\label{dis_er_relation}
\mathcal{L}_{h}^{(i)}E^{(2k+i)}_{i}\leq 0,\ \ i=1,2,\ k = 0,1,2,\dots.
\end{equation}

With the help of (\ref{dis_error_sub}), (\ref{dis_er_relation}) and Theorem \ref{dis_weak_max}, we can prove the convergence of Algorithm \ref{dis_SAM_2} following the same way as that of Theorem \ref{Th:MN_C} and similar convergence results for Algorithm \ref{dis_SAM_2} will be obtained. We just state the results below but omit the details of the proof here.

\begin{Theorem}\label{dis_MN_C}
For a given domain decomposition, if the Schwarz alternating method for 
\begin{equation}\label{dis_L}
\mathcal{L}_{h}W=0\quad\mbox{and}\quad W_{\partial}=0
\end{equation}
is convergent, then the Schwarz alternating algorithm, i.e., Algorithm \ref{dis_SAM_2} for the system (\ref{dis_FOC2}) is convergent. Moreover, if the convergence rate of the Schwarz alternating method for (\ref{dis_L}) is $\rho_{e,d}\in(0,1)$ under the maximum norm, 
then for $k=1, 2, ...$, we have
\begin{equation}\nonumber
\max\{E^{(2k)}\}\leq \rho_{e,d} \max\{E^{(2(k-1))}\}.
\end{equation}
\end{Theorem}

\begin{Remark}
With the help of Theorem \ref{dis_weak_max} and following the way in \cite{PLions_1989}, we can obtain the convergence rate $\rho_{e,d}$ of the Schwarz alternating method for the equation (\ref{dis_L}) under the maximum norm.
\end{Remark}

\begin{Remark}
As in the continuous case, the uniform upper bound $\rho_{e,d}$ here is not optimal, which has been observed in the numerical tests. 
\end{Remark}

Parallel to Corollary \ref{cor:L2_C}, we have the following corollary.
\begin{Corollary} \label{dis_L2_C}
Suppose that the assumptions in Theorem \ref{dis_MN_C} hold, then we have 
\begin{equation}\nonumber
h^2 (E^{(2k)})^T E^{(2k)} = \sum\limits_{s=1}^{(N+1)^2}h^2(E^{(2k)})_{s}\leq C\rho_{e,d}^k\max\{E^{(0)}\},
\end{equation}
where $C>0$ is a constant independent of $h$, $\rho_{e,d}$ and $\alpha$.
\end{Corollary}

\begin{Remark}
The analysis and convergence results in this section are all valid for the central difference scheme in one dimension and the seven-point finite difference scheme in three dimension.
\end{Remark}

\section{Numerical experiments}\label{NE}
In this section we will carry out some numerical experiments to test the performance of the Schwarz alternating method. In our test, we consider the optimal control problem
\begin{eqnarray}\nonumber
\min\limits_{u\in L^2(\Omega)} \
J(y,u)=\frac{1}{2}\|y-y_d\|_{L^{2}(\Omega)}^2 + \frac{\alpha}{2}\|
u\|_{L^{2}(\Omega)}^2
\end{eqnarray}
subject to
\begin{equation}\nonumber
\left\{ \begin{aligned} -\Delta y=f+u \ \ &\mbox{in}\
\Omega, \\
 \ y=0  \ \quad \quad &\mbox{on}\ \partial\Omega,\\
\end{aligned} \right.
\end{equation}
where $\Omega=(0,1)^2$, $y_{d}=\sin(\pi x_1)\sin(\pi x_2)$, $f=2\pi^2 \sin(\pi x_1)\sin(\pi x_2)$. The exact solution of this problem is $y=\sin(\pi x_1)\sin(\pi x_2)$, $p=0$ and $u=0$ for any $\alpha>0$.  For a given uniform grid of $\Omega$ with grid size $h$, we consider the five-point finite difference scheme of this problem which gives the following discrete problem
\begin{eqnarray}\nonumber
\min\limits_{U\in \mathbb{R}^n} \
J(Y,U)=\frac{1}{2}\|Y_{I}-Y_{d,I}\|_{\mathbb{R}^n}^2 + \frac{\alpha}{2}\|
U_{I}\|_{\mathbb{R}^n}^2
\end{eqnarray}
subject to
\begin{equation}\nonumber
\left\{
\begin{aligned}
\mathcal{L}_{h}Y&=F_{I}+U_{I},\\
Y_{\partial}&=0,
\end{aligned}
\right.
\end{equation}
where $Y_{d}$ and $F$ are corresponding to the values of $y_{d}$ and $f$ at the grid points and the meanings of the notations are just as we stated in the previous section. We drop the scaling constant $h^2$ in the objective function in the above formula. By Remark \ref{Remark_1}, this will not affect the performance of the algorithm. We should point out that although we just list the results related to the above discrete optimal control problem, we have also tested the case of random $Y_{d}$, $F$ and it gives similar performances as these results we show here. 

\begin{figure*}[htbp]
\begin{center}
\begin{tikzpicture}[scale = 2.5]
\draw (0,0) -- (2,0);
\draw (0,0) -- (0,2);
\draw (0,2) -- (2,2);
\draw (2,2) -- (2,0);
\foreach \x in {1, 2, ..., 15}
     {
        \draw[color = grey!20] (\x*2/16,0)--(\x*2/16,2);
     }
\foreach \x in {1, 2, ..., 15}
     {
        \draw[color = grey!20] (0, \x*2/16)--(2, \x*2/16);
     }
\draw[dashed] (1,0) -- (1,2);
\draw[dashed, color=red] (6*2/16,0) -- (2,0);
\draw[dashed, color=red] (2,0) -- (2,2);
\draw[dashed, color=red] (6*2/16,0) -- (6*2/16,2);
\draw[dashed, color=red] (2,2) -- (6*2/16,2);
\draw[dashed, color=blue] (10*2/16,0) -- (0,0);
\draw[dashed, color=blue] (0,0) -- (0,2);
\draw[dashed, color=blue] (10*2/16,0) -- (10*2/16,2);
\draw[dashed, color=blue] (10*2/16,2) -- (0,2);
\end{tikzpicture}
\end{center}
\caption{Grid and decomposition of $\Omega$.}
\label{Fig:De_Mesh}
\end{figure*}
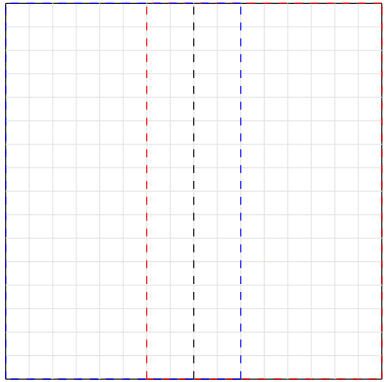

In our test, we decompose the domain in two non-overlapping  parts  first and then we add several layers for each part to obtain  the overlapping domain decomposition (cf. Figure \ref{Fig:De_Mesh}).

In the test, we set the mesh size $h=1/64$ and choose $\alpha=10^{-2},\ 10^{-4},\ 10^{-6}$. We denote the number of the layers added to each subdomain by $\delta$ which can be used to describe the overlapping size of the decomposition of the domain  and choose $\delta=1,\ 2,\ 3,\ 4$, i.e., the overlapping size is $h,\ 2h,\ 3h$ and $4h$ in the test.  We test the impact of the overlapping size and the parameter $\alpha$ on the convergence performance of the algorithm and compare the performances of the algorithm in the equation case and the optimal control problem case.  We get the exact solution of the discretization problem by solving the first order optimality system. $\|E\|_{\infty}$ is the maximum norm of the error between the iterative solution of the Schwarz alternating method and the exact solution of the discrete problem. In the optimal control problem case $\|E\|_\infty := \|E_y^2 + \alpha^{-1}E_p^2\|_\infty$ where $E_y$ and $E_p$ are the error vectors of the state and adjoint respectively. We denote by $\rho_{e,d}$ and $\rho_{c,d}$ the convergence rates of the method for the elliptic equation and the optimal control problem respectively. The test results are reported in Table \ref{table:1}.
  
\begin{figure}
\begin{minipage}{0.48\linewidth}
  \centerline{\includegraphics[width=7.5cm]{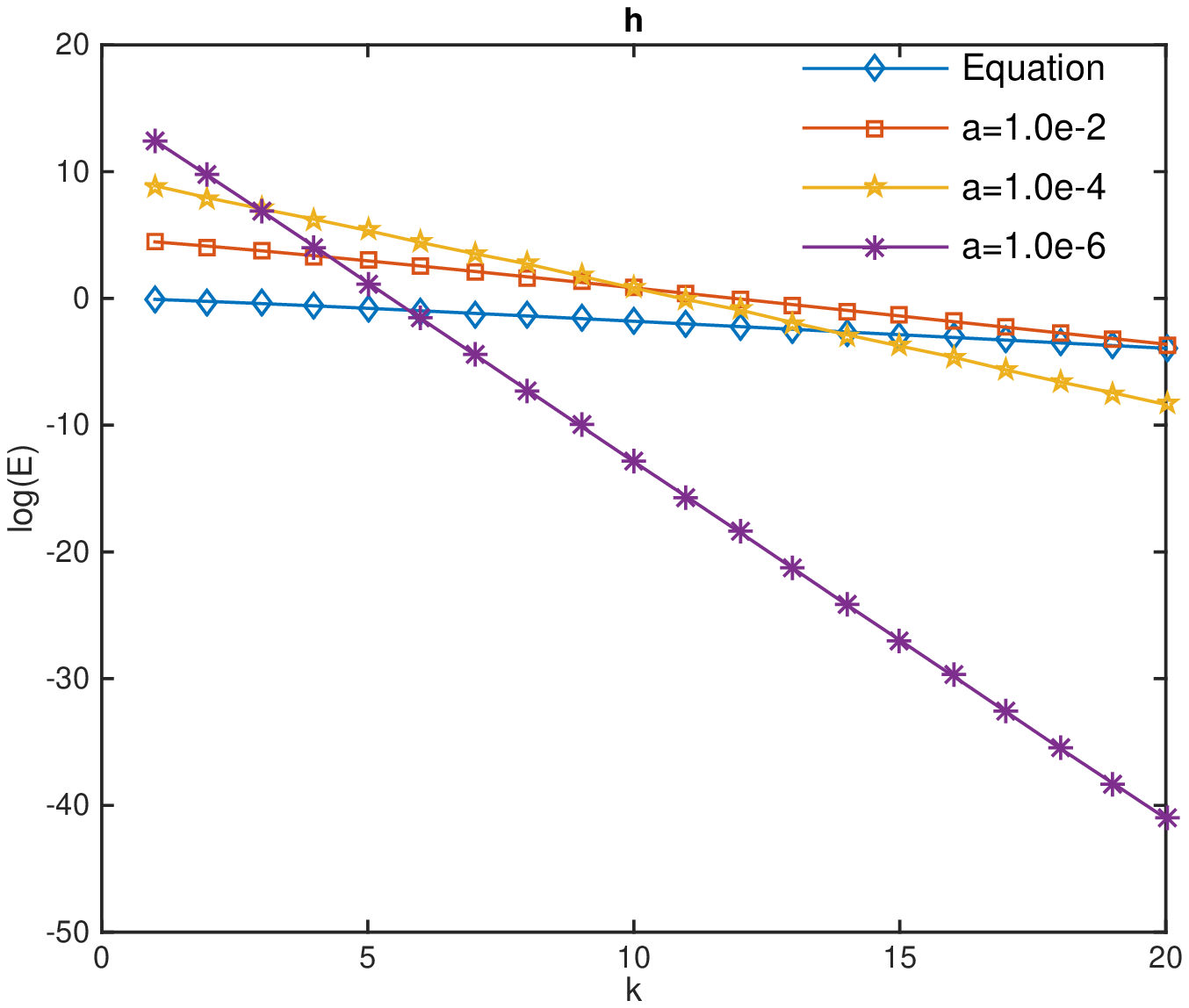}}
  \centerline{(a) $\delta=h$}
\end{minipage}
\hfill
\begin{minipage}{.48\linewidth}
  \centerline{\includegraphics[width=7.5cm]{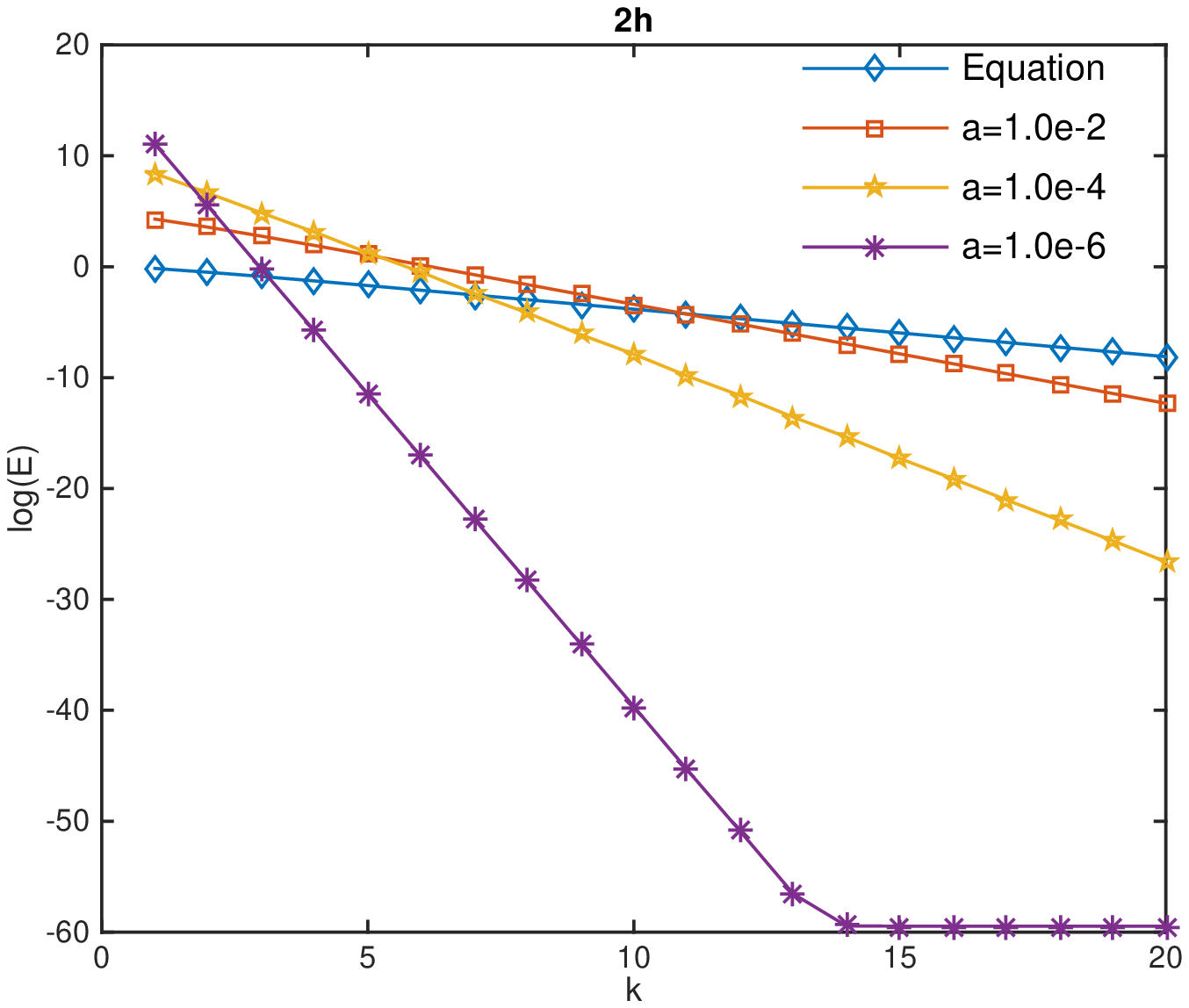}}
  \centerline{(b) $\delta=2h$}
\end{minipage}
\vfill
\begin{minipage}{0.48\linewidth}
  \centerline{\includegraphics[width=7.5cm]{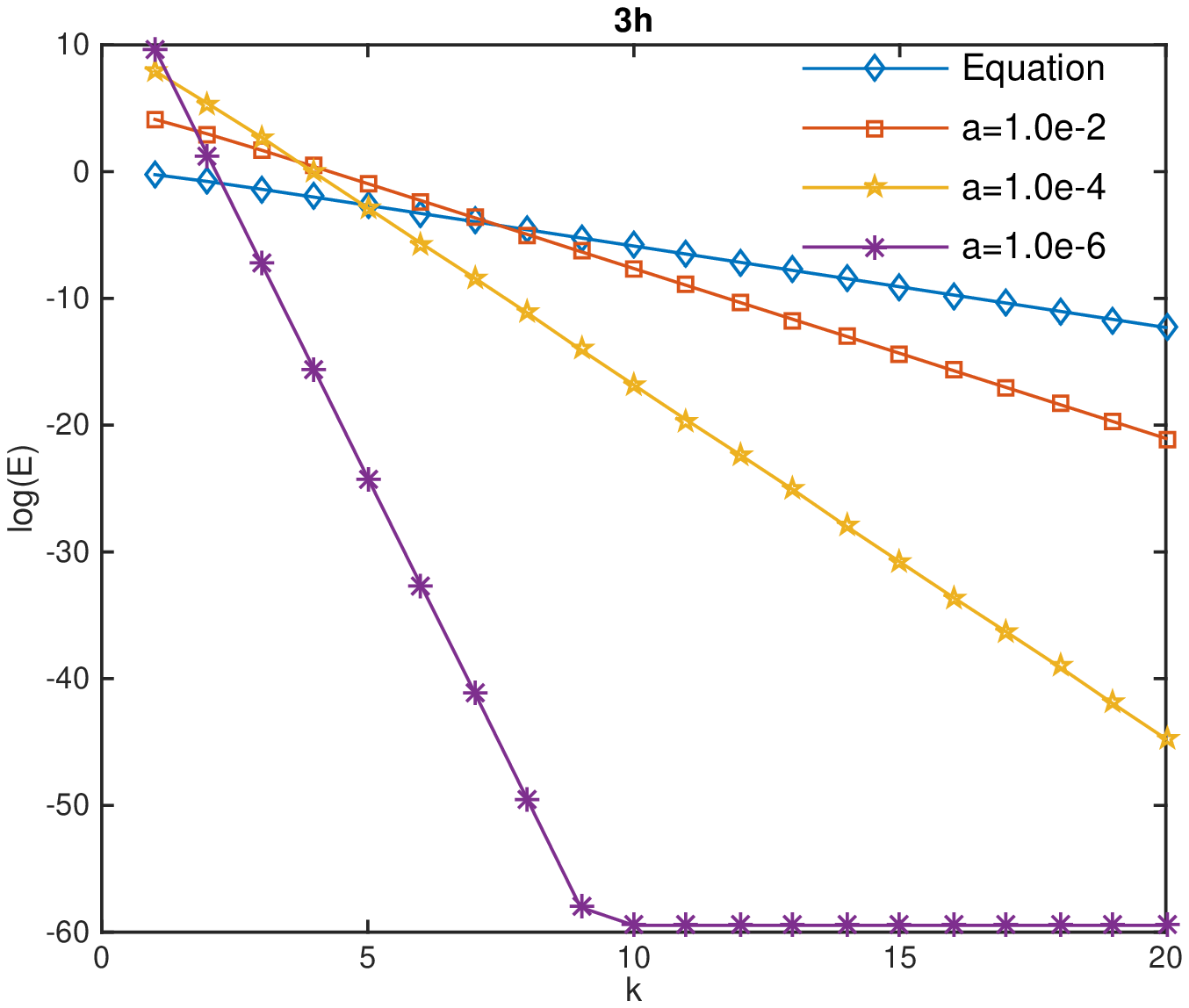}}
  \centerline{(c) $\delta=3h$}
\end{minipage}
\hfill
\begin{minipage}{0.48\linewidth}
  \centerline{\includegraphics[width=7.5cm]{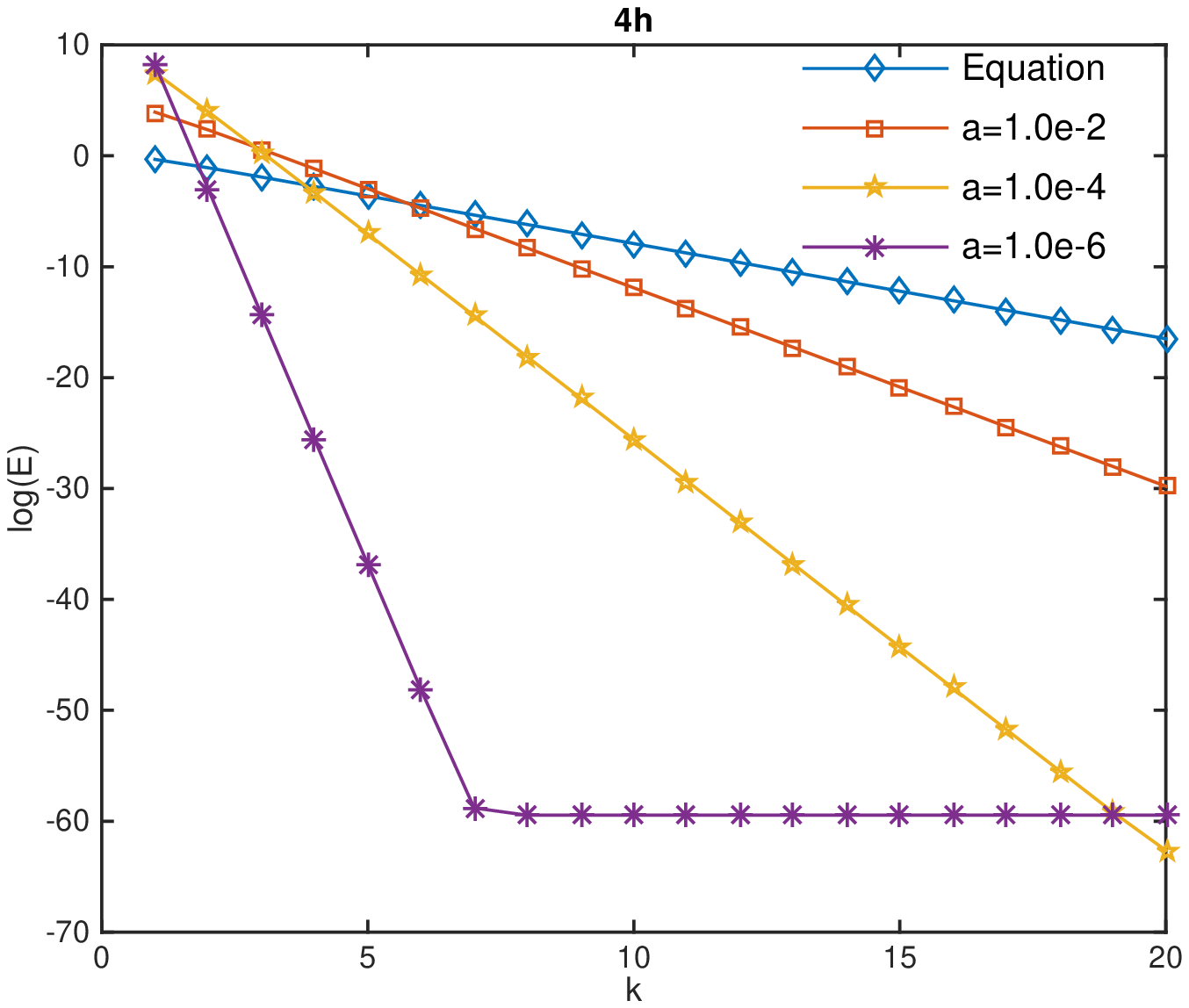}}
  \centerline{(d) $\delta=4h$}
\end{minipage}
\caption{The convergence results}
\label{Figure:con_rate}
\end{figure}

In Table \ref{table:1}, for different overlapping size and different $\alpha$, we present the maximum norm errors of the first five iterations of the Schwarz alternating algorithm  for the elliptic equation and the optimal control problem. We also compute the convergence rate in each test setting. It shows that for the same decomposition of $\Omega$, the convergence rate of the equation case dominates the convergence rate of the optimal control problem case uniformly with respect to $\alpha$. The results also show that $\rho_{c,d}\leq \rho_{e,d}^2$ which has been proved in 1$D$ continuous case (see (\ref{eq:Rate_improve})). In both cases, if the overlap increases, the convergence rate will become smaller.  As $\alpha$ goes smaller, the convergence rate of the method for the optimal control problem will become smaller as well. In order to make the comparison more vivid, we plot the convergence curves of  the algorithm in different test settings and include them in Figure \ref{Figure:con_rate}. These results are in good agreement with our theoretical results. 

\begin{table}
\footnotesize
\centering
\caption{ The errors and convergence rates of the Schwarz alternating method versus $\alpha$, $\delta$ and the iteration number $k$.}\label{table:1}
\begin{tabular}{|*{10}{c|}}
\hline
\multicolumn{2}{|c|}{\ }&\multicolumn{2}{|c|}{Elliptic equation case}&\multicolumn{2}{|c|}{$\alpha=10^{-2}$}&\multicolumn{2}{|c|}{$\alpha=10^{-4}$}&\multicolumn{2}{|c|}{$\alpha=10^{-6}$}\\
\hline
$\delta$ & $k$ & $ \|E\|_{\infty}$ & $\rho_{e,d}$ & $\|E\|_{\infty}$ & $\rho_{c,d}$ & $\|E\|_{\infty}$& $\rho_{c,d}$&$\|E\|_{\infty}$& $\rho_{c,d}$\\
\hline
\multirow{5}{*}{$1$}
&1&9.2738e-1 &- &8.6604e1 &- &7.2285e3 &- &2.5279e5 &-\\
\cline{2-10}
&2&7.9016e-1 &8.5204e-1  &6.1720e1 & 7.1267e-1 &2.7918e3 &3.8623e-1  &1.6005e4 & 6.3314e-2\\
\cline{2-10}
&3&6.6541e-1 &8.4212e-1 &4.2779e1 & 6.9311e-1 &1.1760e3 & 4.2122e-1 &9.8908e2 & 6.1797e-2 \\
\cline{2-10}
&4&5.5466e-1 &8.3357e-1 &2.9006e1 &6.7803e-1 &5.1234e2 &4.3566e-1 &5.9773e1 & 6.0434e-2\\
\cline{2-10}
&5&4.5850e-1 &8.2663e-1 &1.9388e1 &6.6842e-1 &2.1107e2 &4.1198e-1 &3.4274e0 &5.7340e-2\\
\hline
\multirow{5}{*}{$2$}
&1&8.5724e-1 &- & 7.3447e1 &- & 4.5124e3 &- & 6.3505e4 &-\\
\cline{2-10}
&2&6.0803e-1 &7.0929e-1 &3.5317e1 &4.8085e-1 &7.6078e2 &1.6860e-1 & 2.3274e2 & 3.6649e-3 \\
\cline{2-10}
&3&4.1559e-1 &6.8351e-1 &1.5790e1 &4.4708e-1 &1.2923e2 & 1.6986e-1 &8.1402e-1 & 3.4975e-3 \\
\cline{2-10}
&4&2.7760e-1 &6.6795e-1 &6.7994e0 &4.3063e-1 &2.2636e1 &1.7516e-1 &3.0695e-3 &3.7709e-3 \\
\cline{2-10}
&5&1.8311e-1 &6.5963e-1 &2.8540e0 &4.1974e-1 &3.5541e0 & 1.5701e-1 &1.0338e-5 & 3.3680e-3\\
\hline
\multirow{5}{*}{$3$}
&1&7.8980e-1 &- &6.1795e1 &- &2.7926e3 &- &1.6005e4 &-\\
\cline{2-10}
&2&4.5793e-1 &5.7981e-1 &1.9350e1 &3.1314e-1 &2.1125e2 &7.5645e-2 &3.4274e0 &2.1414e-4 \\
\cline{2-10}
&3&2.5008e-1 &5.4609e-1 &5.4692e0 &2.8264e-1 &1.5174e1 &7.1830e-2 &7.2018e-4 &2.1012e-4 \\
\cline{2-10}
&4&1.3321e-1 &5.3268e-1 &1.4701e0 &2.6879e-1 &9.0707e-1 &5.9778e-2 &1.5694e-7 &2.1792e-4 \\
\cline{2-10}
&5&7.0335e-2 &5.2800e-1 &3.8813e-1&2.6402e-1 &5.7030e-2 &6.2872e-2 &3.1882e-11 &2.0315e-4 \\
\hline
\multirow{5}{*}{$4$}
&1&7.2530e-1 &- &5.1592e1 &- &1.8031e3 &- &3.8839e3 &-\\
\cline{2-10}
&2&3.3954e-1 &4.6814e-1 &1.0356e1 &2.0074e-1 &5.3161e1 &2.9483e-2 &5.0780e-2 &1.3074e-5 \\
\cline{2-10}
&3&1.4774e-1 &4.3511e-1 &1.8249e0 &1.7621e-1 &1.4308e0 &2.6915e-2 &6.2383e-7 & 1.2285e-5\\
\cline{2-10}
&4&6.3005e-2 &4.2647e-1 &3.0868e-1 &1.6915e-1 &3.5916e-2 &2.5101e-2 &7.9164e-12 &1.2690e-5 \\
\cline{2-10}
&5&2.6743e-2 &4.2446e-1 &5.1728e-2 &1.6758e-1 &8.7228e-4 &2.4287e-2 &1.0084e-16 &1.2739e-5 \\
\hline
\end{tabular}
\end{table}

We also compute the convergence results of the Schwarz alternating method for the $\alpha$-dependent elliptic equation
\begin{equation}\label{eq:alpha_elliptic}
-\Delta \tilde{w} + 2\alpha^{-\frac{1}{2}}\tilde{w} = 0\ \ \mbox{in}\ \ \Omega\quad \mbox{and}\quad \tilde{w} = 0\ \ \mbox{on}\ \ \partial\Omega
\end{equation}
under the same settings as before. The numerical results are reported in Table  \ref{table:better_estimate}. The results for the case where $\alpha = 10^{-6}$ in Table \ref{table:1} and Table \ref{table:better_estimate} show that the convergence rate of the Schwarz alternating method for the $\alpha$-dependent equation is quite a good estimate of the convergence rate of the method for  the optimal control problem when $\alpha$ is small, which agrees with the estimate in (\ref{eq:better_estimate}).

\begin{table}
\footnotesize
\centering
\caption{ The errors and convergence rates of the Schwarz alternating method versus $\delta$ and the iteration number $k$ for the $\alpha$-dependent elliptic equation with $\alpha = 10^{-6}$.}\label{table:better_estimate}
\begin{tabular}{|*{9}{c|}}
\hline
\multirow{2}{*}{\diagbox{$k$}{$\delta$}} &
\multicolumn{2}{|c|}{$1$}&\multicolumn{2}{|c|}{$2$}&\multicolumn{2}{|c|}{$3$}&\multicolumn{2}{|c|}{$4$}\\
\cline{2-9}
& $ \|E\|_{\infty}$ & $\rho_{e,d}$ & $\|E\|_{\infty}$ & $\rho_{e,d}$ & $\|E\|_{\infty}$& $\rho_{e,d}$&$\|E\|_{\infty}$& $\rho_{e,d}$\\
\hline
1& 2.5396e-1 &- & 6.4497e-2 &- & 1.6380e-2 &- & 4.1599e-3 &-\\
\hline
2& 1.6380e-2 & 6.4497e-2 &2.6830e-4 &4.1599e-3 & 4.3948e-6 &2.6830e-4 & 7.1986e-8 & 1.7305e-5\\
\hline
3& 1.0565e-3 & 6.4497e-2 &1.1161e-6 &4.1599e-3 & 1.1791e-9 &2.6830e-4 & 1.2456e-12 & 1.7303e-5\\
\hline
4& 6.8139e-5 & 6.4497e-2 &4.6428e-9 &4.1599e-3 & 3.1632e-13 &2.6827e-4 & 2.1539e-17 & 1.7293e-5\\
\hline
5& 4.3948e-6 & 6.4497e-2 & 1.9313e-11&4.1599e-3 & 8.4825e-17 &2.6816e-4 & 3.7189e-22 & 1.7266e-5\\
\hline
\end{tabular}
\end{table}

\section{Conclusion}
By using the maximum principle of the elliptic operator, we can analyze the Schwarz alternating method for unconstrained elliptic optimal control problems and prove the convergence of the method under the maximum norm with a convergence rate which is uniformly bounded by a constant less than one. This will also give a first glimpse of exploring the convergence properties of other overlapping/non-overlapping domain decomposition methods for optimal control problems in both continuous and discrete settings, such as the one-level/two-level additive Schwarz method and the multiplicative Schwarz method and so on, which are originated from this method.

As mentioned in Remark \ref{rem:Th_rate}, the estimated uniform upper bound of the convergence rate in this paper is not optimal. How to improve it to the optimal one (see (\ref{eq:Rate_improve})) and how to give a fine estimate of it (see (\ref{eq:better_estimate})) are very interesting and a little bit tricky problems. Besides, 
the approach in this paper can possibly be extended to local optimal control problems where the control only acts on a subset of the domain rather than the whole domain and the problem with constraints. It is also interesting and important to prove the geometric convergence property of the method for the finite element discretization case, which has been observed numerically (cf. \cite{Tan_phd}). These problems are among our ongoing projects.

\begin{appendix}
\section{The convergence results in one dimensional case}\label{appendix:A}
\setcounter{equation}{0}
In this part, we give an explicit formulation of the convergence rate of the Schwarz alternating method in the one dimensional case. This part is mainly taken from \cite{Tan_phd}. 

We consider the one dimensional 
unconstrained elliptic  distributed optimal  control problem:
\begin{eqnarray}\label{OPTC_1}
\min\limits_{u\in L^2(\Omega)} \
J(y,u)=\frac{1}{2}\|y-y_d\|_{L^{2}(\Omega)}^2 + \frac{\alpha}{2}\|
u\|_{L^{2}(\Omega)}^2
\end{eqnarray}
subject to
\begin{equation}\label{OPT_PDE_1}
\left\{ \begin{aligned} -y''=f+u \ \ &\mbox{in}\
\Omega, \\
  y=0  \ \quad \quad &\mbox{on}\ \partial\Omega,\\
\end{aligned} \right.
\end{equation}
where $\Omega=(0,1)$,
$u\in L^{2}(\Omega)$ is the control variable, $y_d\in L^2(\Omega)$ is the
desired state or observation and $\alpha>0$ is the regularization
parameter.  

We set $\Omega_{1}=(0,s)$, $\Omega_{2}=(r,1)$ with $0<r<s<1$. Then the error systems (\ref{error_system}) are given by
\begin{equation}\label{error_system_1}
\left\{\begin{aligned}
&-(e_{y}^{(2k+1)})'' =-\alpha^{-1} e_{p}^{(2k+1)} \ &\mbox{in}\
(0,s),\  e_{y}^{(2k+1)}(0)=0,\ e_{y}^{(2k+1)}(s)=e_{y}^{(2k)}(s),\\
&-(e_{p}^{(2k+1)})'' = e_{y}^{(2k+1)} \ &\mbox{in}\
(0,s),\  e_{p}^{(2k+1)}(0)=0,\ e_{p}^{(2k+1)}(s)=e_{p}^{(2k)}(s)\\
\end{aligned}
\right.  
\end{equation}
and
\begin{equation}\label{error_system_2}
\left\{\begin{aligned}
&-(e_{y}^{(2k+2)})'' =-\alpha^{-1} e_{p}^{(2k+2)} \ &\mbox{in}\
(r,1),\ e_{y}^{(2k+2)}(r)=e_{y}^{(2k+1)}(s),\  e_{y}^{(2k+2)}(1)=0,\\
&-(e_{p}^{(2k+2)})'' = e_{y}^{(2k+2)} \ &\mbox{in}\
(r,1),\ e_{p}^{(2k+2)}(s)=e_{p}^{(2k+1)}(s),\  e_{p}^{(2k+2)}(1)=0.\\
\end{aligned}
\right.  
\end{equation}

We can solve the systems (\ref{error_system_1}) and (\ref{error_system_2}) directly and obtain the analytic solutions.  Let $\gamma=\frac{\sqrt{2}}{2}\alpha^{-\frac{1}{4}}$ and
\begin{equation}\nonumber
\begin{aligned}
q_{1}&=\sinh(\gamma(s-x))\sin(\gamma(s+x))-\sin(\gamma(s-x))\sinh(\gamma(s+x)),\\
q_{2}&=\cosh(\gamma(s-x))\cos(\gamma(s+x))-\cos(\gamma(s-x))\cosh(\gamma(s+x)),\\
q_{3}&=\cosh(2\gamma s)-\cos(2\gamma s)=2(\cosh^2(\gamma s)-\cos^2(\gamma s))=2(\sinh^2(\gamma s)+\sin^2(\gamma s)),\\
\tilde{q}_{1}&=\sinh(\gamma((1-r)-(1-x)))\sin(\gamma((1-r)+(1-x)))\nonumber\\
&-\sin(\gamma((1-r)-(1-x)))\sinh(\gamma((1-r)+(1-x))),\\
\tilde{q}_{2}&=\cosh(\gamma((1-r)-(1-x)))\cos(\gamma((1-r)+(1-x)))\nonumber\\
&-\cos(\gamma((1-r)-(1-x)))\cosh(\gamma((1-r)+(1-x))),\\
\tilde{q}_{3}&=\cosh(2\gamma (1-r))-\cos(2\gamma (1-r))=2(\sinh^2(\gamma (1-r))+\sin^2(\gamma (1-r))),\\
\end{aligned}
\end{equation}
where $\sinh(x)$, $\cosh(x)$ are the hyperbolic sine, cosine functions 
\begin{equation}\nonumber
\sinh(x)=\frac{e^{x}-e^{-x}}{2}\quad\mbox{and}\quad \cosh(x)=\frac{e^{x}+e^{-x}}{2}.
\end{equation}

The solutions of (\ref{error_system_1}) and (\ref{error_system_2}) are
\begin{equation}\nonumber
e_{y}^{(2k+1)}(x)=-\frac{q_{2}e_{y}^{(2k)}(s)-2\gamma^2 q_{1}e_{p}^{(2k)}(s)}{q_{3}},
\end{equation}
\begin{equation}\nonumber
e_{p}^{(2k+1)}=-\frac{q_{1}e_{y}^{(2k)}(s)+2\gamma^2 q_{2}e_{p}^{(2k)}(s)}{2\gamma^2q_{3}},
\end{equation}
\begin{equation}\nonumber
e_{y}^{(2k+2)}(x)=-\frac{\tilde{q}_{2}e_{y}^{(2k+1)}(r)-2\gamma^2 \tilde{q}_{1}e_{p}^{(2k+1)}(r)}{q_{3}},
\end{equation}
\begin{equation}\nonumber
e_{p}^{(2k+2)}(x)=-\frac{\tilde{q}_{1}e_{y}^{(2k+1)}(r)+2\gamma^2 \tilde{q}_{2}e_{p}^{(2k+1)}(r)}{2\gamma^2\tilde{q}_{3}}.
\end{equation}

Then we have
\begin{equation}\nonumber
\left\{
\begin{aligned}
&(e_{y}^{(2k+1)}(x))^2+\alpha^{-1}(e_{p}^{(2k+1)}(x))^2=L(x,s)\left((e_{y}^{(2k)}(s))^2+\alpha^{-1}(e_{p}^{(2k)}(s))^2\right)\quad x\in [0,s],\\
&(e_{y}^{(2k+2)}(x))^2+\alpha^{-1}(e_{p}^{(2k+2)}(x))^2=R(r,x)\left((e_{y}^{(2k+1)}(r))^2+\alpha^{-1}(e_{p}^{(2k+1)}(r))^2\right)\quad x\in [r,1],\\
\end{aligned}
\right.
\end{equation}
where
\begin{equation}
\left\{
\begin{aligned}
&L(x,s)=\frac{\sinh^2(\gamma x)+\sin^2(\gamma x)}{\sinh^2(\gamma s)+\sin^2(\gamma s)}\quad \ x\in [0,s],\\
&R(r,x)=\frac{\sinh^2(\gamma (1-x))+\sin^2(\gamma (1-x))}{\sinh^2(\gamma (1-r))+\sin^2(\gamma (1-r))}\quad x\in [r,1].
\end{aligned}
\right.
\end{equation}
Furthermore, we have
\begin{equation}\nonumber
\begin{aligned}
&(e_{y}^{(2k+1)}(r))^2+\alpha^{-1}(e_{p}^{(2k+1)}(r))^2=L(r,s)\left((e_{y}^{(2k)}(s))^2+\alpha^{-1}(e_{p}^{(2k)}(s))^2\right)\\
\end{aligned}
\end{equation}
and
\begin{equation}\nonumber
\begin{aligned}
&(e_{y}^{(2k+2)}(s))^2+\alpha^{-1}(e_{p}^{(2k+2)}(s))^2=R(r,s)\left((e_{y}^{(2k+1)}(r))^2+\alpha^{-1}(e_{p}^{(2k+1)}(r))^2\right).\\
\end{aligned}
\end{equation}

Since 
\begin{equation}\nonumber
g(z)=\sinh^2(z)+\sin^2(z)=\frac{1}{2}(\cosh(2z)-\cos(2z))
\end{equation}
is positive, strictly convex, strictly monotone increasing in $(0,\infty)$  and $\lim\limits_{z\rightarrow +\infty}g(z)=+\infty$, we have
\begin{equation}\nonumber
\begin{aligned}
&\max\limits_{x\in [0,1]}\left((e_{y}^{(2k+1)}(x))^2+\alpha^{-1}(e_{p}^{(2k+1)}(x))^2\right)=(e_{y}^{(2k)}(s))^2+\alpha^{-1}(e_{p}^{(2k)}(s))^2,\\
&\max\limits_{x\in [0,1]}\left((e_{y}^{(2k+2)}(x))^2+\alpha^{-1}(e_{p}^{(2k+2)}(x))^2\right)=(e_{y}^{(2k+1)}(r))^2+\alpha^{-1}(e_{p}^{(2k+1)}(r))^2.\\
\end{aligned}
\end{equation}

Hence
\begin{equation}\nonumber
\begin{aligned}
&\max\limits_{x\in [0,1]}\left((e_{y}^{(2k+2)}(x))^2+\alpha^{-1}(e_{p}^{(2k+2)}(x))^2\right)=L(r,s)R(r,s)\max\limits_{x\in [0,1]}\left((e_{y}^{(2k)}(x))^2+\alpha^{-1}(e_{p}^{(2k)}(x))^2\right),\\
\end{aligned}
\end{equation}
where $k=0,1,2,...$.
This means the convergence rate $\rho_{c}=L(r,s)R(r,s)$ and $\rho_c<1$.

It is easy to derive that the convergence rate of the equation case is 
\begin{equation}\nonumber
\rho_{e}=\frac{r(1-s)}{s(1-r)}.
\end{equation}

By the convexity of $g(z)$, $g(0)=0$ and the fact $0<r<s$, we have 
\begin{equation}\nonumber
\cosh(2\gamma r)-\cos(2\gamma r)= \cosh(2\gamma s\frac{r}{s} )-\cos(2\gamma s \frac{r}{s} )<\frac{r}{s}(\cosh(2\gamma s)-\cos(2\gamma s)).
\end{equation}
It implies 
$$\rho_{c}<\rho_{e}.$$

\begin{figure}
\centerline{\includegraphics[width=7.5cm]{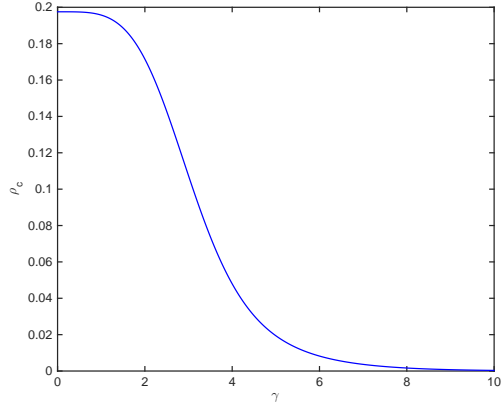}}
\caption{Impact of $\gamma = \frac{\sqrt{2}}{2}\alpha^{-\frac{1}{4}}$ with $r=0.4$, $s=0.6$.}
\label{Fig:Rate_gamma}
\end{figure}

Furthermore, we can prove that for given $0<r<s<1$, $\rho_c$ is strictly monotone decreasing with respect to $\gamma$ (cf. Figure \ref{Fig:Rate_gamma} for an illustration).  Since $\gamma=\frac{\sqrt{2}}{2}\alpha^{-\frac{1}{4}}$, we know that as $\alpha$ decreasing, $\rho_c$ will decrease for fixed $0<r<s<1$. Our numerical results for the two dimensional case also corroborate this observation.

Note that
\begin{equation}\nonumber
\frac{g(\gamma r)}{g(\gamma s)}\approx e^{2\gamma (r-s)}\ \ \mbox{as}\ \ \gamma \rightarrow +\infty\quad \mbox{and}\quad \lim\limits_{\gamma\rightarrow 0^{+}}\frac{g(\gamma r)}{g(\gamma s)} = \frac{r^2}{s^2}.
\end{equation}
It follows
\begin{equation}\label{eq:asym_rho}
\rho_c\approx e^{2\sqrt{2}\alpha^{-\frac{1}{4}} (r-s)}\ \ \mbox{as}\ \ \alpha\rightarrow 0^{+}\quad \mbox{and}\quad \rho_c \approx \rho_e^2\ \ \mbox{as}\ \ \alpha\rightarrow +\infty,
\end{equation}
which give the asymptotic behaviors of the convergence rate as $\alpha$ goes to zero and infinity. Furthermore, by applying the monotonicity of $\rho_c$, we have for any $\alpha >0$
\begin{equation}\label{eq:Rate_improve}
\rho_c \leq \rho_e^2.
\end{equation}

\section{A better estimate in one dimensional case}
\setcounter{equation}{0}
In Appendix \ref{appendix:A}, we compare the convergence rate of the Schwarz alternating method for the elliptic equation
$$ -w'' =0\ \ \mbox{in}\ (0,1)\quad \mbox{and}\quad \ w(0) = w(1) =0$$
and the convergence rate of the method for the optimal control problem. Now we consider the convergence rate of the method for the elliptic equation
$$ -\tilde{w}'' + \beta^2 \tilde{w} =0\ \ \mbox{in}\ (0,1)\quad \mbox{and}\quad \ \tilde{w}(0) = \tilde{w}(1) =0,$$
where $\beta\geq 0$.

We use the same setting as those in Appendix \ref{appendix:A}. In this case, a direct calculation gives the convergence rate of the method as
\begin{equation}\nonumber
\tilde{\rho}_{e,\beta} = \frac{\sinh(\beta r)\sinh(\beta (1-s))}{\sinh(\beta s)\sinh(\beta (1-r))}.
\end{equation}

Since $\sinh(\beta x)$ is positive, strictly convex, strictly monotone increasing in $(0,\infty)$ and $\sinh(0) = 0$, 
we have
\begin{equation}\nonumber
 \sinh(\beta r) = \sinh(\beta s\cdot\frac{r}{s}) < \frac{r}{s} \sinh(\beta s)
\end{equation}
which implies
\begin{equation}\nonumber
 \frac{\sinh(\beta r)}{\sinh(\beta s)} < \frac{r}{s}
\end{equation}
and
$$ \tilde{\rho}_{e,\beta}< \rho_e\quad \forall\ \beta >0.$$

Note that
$$ \tilde{\rho}_{e,\beta} \approx e^{2\beta(r-s)}\quad \mbox{as}\ \ \beta\rightarrow +\infty\quad \mbox{and}\quad \lim\limits_{\beta\rightarrow 0^+}\tilde{\rho}_{e,\beta} = \rho_e.$$

If we take $\beta = 2\gamma$, i.e., $\beta = \sqrt{2}\alpha^{-\frac{1}{4}}$, we have
\begin{equation}\label{eq:better_estimate}
\rho_c \approx \tilde{\rho}_{e, \sqrt{2}\alpha^{-\frac{1}{4}}}\quad \mbox{as}\ \ \alpha\rightarrow 0^+\quad \mbox{and}\quad \rho_c \approx \tilde{\rho}_{e, \sqrt{2}\alpha^{-\frac{1}{4}}}^2\ \ \mbox{as}\ \ \alpha\rightarrow +\infty.
\end{equation}


\end{appendix}


 \medskip

\end{document}